\documentclass[final,leqno]{elsarticle}

\usepackage{amsmath,amssymb}
\usepackage{algpseudocode}
\usepackage{graphicx}% Include figure files
\usepackage{dcolumn}% Align table columns on decimal point
\usepackage{booktabs}
\usepackage{bm}% bold math
\usepackage[font=footnotesize,format=plain,labelfont=bf]{caption}
\usepackage{subcaption}
\usepackage{xspace}
\usepackage{wrapfig}
\usepackage{tikz}
\usepackage{longtable}
\usepackage{multirow}
\usepackage{hyperref}

 \newtheorem{remark}{Remark}
 \newtheorem{definition}{Definition}
\usepackage{float}

\newcommand{\pvct}[1]{\bm{#1}}
\newcommand{\vct}[1]{\bm{\mathsf{#1}}}

\newcommand{\pxx}{\pvct{x}}

\newcommand{\uu}{\vct{u}}

\newcommand{\mtx}[1]{\bm{\mathsf{#1}}}

\newcommand{\mtwo}[4]{\left[\begin{array}{cc} #1 & #2 \\ #3 & #4 \end{array}\right]}
\newcommand{\vtwo}[2]{\left[\begin{array}{cc} #1 \\ #2 \end{array}\right]}

\newcommand{\pgnotate}[1]{}

\newcommand{\IncImpHmgSolnOp}{\mtx{\Phi}}
\newcommand{\OutImpPrtSolnOp}{\mtx{\Gamma}}

\newcommand{\HmgSolnOp}{\mtx{\Psi}}
\newcommand{\PrtSolnOp}{\mtx{Y}}

\newcommand{\outimpprtvec}{\vct{h}}
\newcommand{\prtsoln}{\tilde{u}}
\newcommand{\hmgsoln}{\bar{u}}
\newcommand{\prtsolnvec}{\tilde{\vct{u}}}
\newcommand{\hmgsolnvec}{\vct{\bar{u}}}

\definecolor{lightblue}{rgb}{0.39, 0.58, 0.93}
\definecolor{darkmagenta}{rgb}{0.55, 0.0, 0.55}
\definecolor{green}{rgb}{0.24, 0.71, 0.54}
\definecolor{magenta}{rgb}{0.73, 0.33, 0.83}
\definecolor{darkblue}{rgb}{0.1, 0.1, 0.44}

\title{A parallel shared-memory implementation of a high-order accurate solution technique for variable coefficient Helmholtz problems}

\author[rice]{Natalie N.~Beams}
\ead{natalie.beams@rice.edu}

\author[rice]{Adrianna Gillman\corref{cor1}}
\ead{adrianna.gillman@rice.edu}

\author[total,VTech]{Russell J. Hewett}
\cortext[cor1]{Corresponding author}
\address[rice]{Department of Computational and Applied Mathematics, Rice University}
\address[total]{Total E\&P Research \& Technology USA  }
\address[VTech]{Department of Mathematics, Virginia Tech}
\begin{document}
 \begin{frontmatter}

\begin{abstract} 
The recently developed Hierarchical Poincar\'e-Steklov (HPS) method 
is a high-order discretization technique that comes with a direct
solver.  Results from previous papers demonstrate the method's ability 
 to solve Helmholtz problems to high 
accuracy without the so-called pollution effect.  While the asymptotic 
scaling of the direct solver's computational cost is the same as the
nested dissection method, serial implementations of the solution technique
are not practical for large scale numerical simulations.  This 
manuscript presents the first parallel implementation of the 
HPS method.  Specifically, we introduce an approach for a shared memory 
implementation of the solution technique utilizing parallel linear 
algebra.  This approach is the foundation for future
large scale simulations on supercomputers and clusters
with large memory nodes.  Performance results on a 
desktop computer (resembling a large memory node) are 
presented.

\end{abstract}

\begin{keyword}Numerical partial differential equations \sep Direct solver \sep 
High-order discretization \sep Nested dissection \sep OpenMP \sep shared-memory parallelization \sep MKL\end{keyword}

%\begin{AMS} TODO -- ams keywords    \end{AMS}

\end{frontmatter}

\section{Introduction}
Consider the variable coefficient Helmholtz problem 
\begin{equation}
\left\{\begin{aligned}
-\Delta u(\vct{x}) - \kappa^2 c(\vct{x})u(\vct{x})=&\ s(\vct{x})\qquad &\vct{x} \in \Omega,\\
\frac{\partial u}{\partial \vct{\nu}} + \textrm{i}\eta u    =&\ t(\vct{x})\qquad &\vct{x} \in \Gamma = \partial \Omega,
\end{aligned}\right.
\label{eq:basic}
\end{equation}
where $\Omega$ is a rectangle in $\mathbb{R}^2$,
$\kappa$ is the wave number, $\vct{\nu}$ is the outward facing normal on $\Gamma$,
$\eta\in \mathbb{C}$ ($\Re(\eta)\neq 0$) and $u(\vct{x})$ is the unknown solution.  The functions 
$s(\vct{x})$, $t(\vct{x})$, and $c(\vct{x})$ are assumed to be smooth.  
Solutions to this boundary value problem are oscillatory and the frequency at which the solutions oscillate is dictated by $\kappa$.  In other words, as $\kappa$ grows, the solution becomes more oscillatory.  The task of creating high-order approximate solutions to boundary value problems of the form (\ref{eq:basic}),  
where the number of discretization points per wavelength is fixed has been a 
challenge for some time.  The recently developed Hierarchical Poincar\'e-Steklov (HPS) method 
is a high-order discretization technique that comes with an efficient direct solver and does not, in numerical experiments, suffer from the so-called \textit{pollution} effect
\cite{2013_martinsson_ItI}.  For the HPS method to be useful for large scale 
computations, a high performance computing implementation of the method is necessary.  This paper presents the first such implementation.  The 
implementation is for a shared memory machine that is representative of the large
memory nodes in upcoming supercomputers and clusters.

While this paper considers the Helmholtz impedance boundary value problem (\ref{eq:basic}) for
simplicity of presentation, the technique can be used to solve problems with 
other boundary conditions with minor modifications (see \cite{2013_martinsson_ItI}). 
Additionally, the parallelization technique can be 
applied directly to the variant of the HPS method for elliptic boundary value problems \cite{2012_spectralcomposite, 2016_bodyload}.

\subsection{Overview of discretization technique}
\label{sec:overview}
Roughly speaking, the discretization technique and construction of the direct solver 
can be broken into three steps:\\

\begin{tabular}{ll}
Step 1:& The geometry is partitioned into a collection 
 of disjoint patches\\
 & sized so that a boundary value problem on the patch can be \\
 & solved to high accuracy via a classic spectral collocation method \\
 & (e.g. \cite{2000_trefethen_spectral_matlab}).\\
Step 2:& Each patch is discretized using a high order spectral collocation\\ 
 & technique.  Approximate boundary (Poincar\'e-Steklov) and \\
 & solution operators are constructed.  \\
 Step 3:& In a hierarchical fashion, the patches are ``glued'' 
 together two at  \\
 & a time by enforcing continuity conditions on the solution \\
 &via the Poincar\'e-Steklov operators on the boundaries
 of each  \\ 
&patch.  For each merged patch, corresponding boundary and \\
&solution operators are constructed. 
\end{tabular}

These three steps comprise the \textit{precomputation} stage of the solution 
technique.  Once the precomputation is complete, the task of finding the solution to (\ref{eq:basic}) for 
any choice of body load $s(\vct{x})$ and boundary data $t(\vct{x})$ 
is simply a collection of small matrix vector multiplies involving the precomputed
operators. This is called the \textit{solve} stage.

The domain decomposing nature of the algorithm provides significant 
opportunities for parallelism.  For two-dimensional problems, the required 
dense linear algebra involves matrices corresponding to operators that live on 
a line, which keeps the overall computational cost in FLOPs low. The distribution of the work while moving through the 
hierarchical tree in both stages of the algorithm are explored in this paper.

While the method can be employed with any Poincar\'e-Steklov operator, this paper 
uses the impedance-to-impedance (ItI) operator, which is ideal for Helmholtz problems. For general elliptic problems, 
the Dirichlet-to-Neumann operator is a suitable choice 
\cite{2013_martinsson_DtN_linearcomplexity, 2012_spectralcomposite, 2016_bodyload}.

\subsection{Related to prior work}
The original HPS method \cite{2012_spectralcomposite} was designed for
elliptic partial differential equations and the local discretization 
utilized classic spectral collocation techniques, which involved 
points at the corners of leaf boxes.  These corner discretization 
points were not ideal for the ``gluing" procedure in Step 2.  In \cite{2013_martinsson_DtN_linearcomplexity,
2016_bodyload,2013_martinsson_ItI} corner points were removed by using interpolation 
operators to represent the boundary operators only on edges of boxes.  Most 
recently, in \cite{2018_geldermans_gillman_adaptive}, a new spectral collocation 
scheme is presented which does not place any discretization points the corners
of boxes.  The parallelization of this latest version of the method is presented
in this paper.

% The development of the HPS method has evolved from using a classic spectral collocation 
% scheme on the small patches to the current version which has a special spectral basis 
% specifically for the HPS method which does not have corner points.  This work uses the 
% new spectral collocation technique presented in \cite{2018_geldermans_gillman_adaptive}.  

The direct solver for the HPS discretization is related to 
the direct solvers for sparse systems arising from finite 
difference and finite element discretizations of elliptic
PDEs, such as the classical nested dissection method of 
George \cite{george_1973} and the multifrontal 
methods by Duff et al. \cite{1989_directbook_duff}.  These 
methods can be viewed as a hierarchical version of the
``static condensation'' idea in finite element analysis 
\cite{1974_wilson_static_condensation}.  High-order finite 
difference and finite element discretizations lead to large 
frontal matrices, and consequently very high cost of the 
LU-factorization (see, e.g., Table 2 in \cite{2013_martinsson_DtN_linearcomplexity}).  
It has been demonstrated that the dense matrices that arise in these solvers
have internal structure that allows the direct solver to be accelerated to linear or close to linear
complexity, see, e.g., \cite{2009_xia_superfast,Adiss,2007_leborne_HLU,2009_martinsson_FEM,2011_ying_nested_dissection_2D}.
The two-dimensional HPS solution technique has one dimensional ``dividers'' independent of order and thus
the direct solver only pays (in terms of computational complexity) the price of the high-order discretization 
at the lowest level in the hierarchical tree.  
The same ideas that accelerate the nested dissection and multifrontal solvers can be applied 
Helmholtz problems, the scaling of these accelerated solvers deteriorates. It 
should be noted that the parallelization technique presented in this 
manuscript does not apply to the linear scaling direct solver version in \cite{2013_martinsson_DtN_linearcomplexity}.

There are multiple widely-available libraries for high performance 
parallel implementations of direct factorization for sparse matrices, i.e.~nested dissection and its variants. 
SuperLU \cite{superlu_ug99, superlu_mt-1999, lidemmel-superlu-dist-2003}
takes either
a left-looking (shared memory) or right-looking (distributed memory) approach to 
factorization.  To minimize idle time when the number of independent                  
tasks is less than the number of available processors, SuperLU implements pipelining,
where portions of dependent tasks are computed simultaneously and waiting only occurs when a task
cannot continue without receiving necessary data from a related task.
The multifrontal solvers in WSMP \cite{gupta-wsmp-sparse-solver-2007} and MUMPS \cite{mumps-2001,mumps-dynamic-scheduling-2006} both employ multiple strategies for parallelism
based on the hierarchical nature of the multifrontal algorithm.  First, there is parallelism from 
the natural independent
calculation of subtrees that do not depend on each other.  As the elimination continues and 
the number of independent calculations is greater than the number of processes,         
the processes share the calculations through parallel linear algebra. This split between types of parallelism is also employed for HSS matrices in STRUMPACK \cite{Strumpack}.  The tailoring of parallelism 
to the algorithm's tree structure is, in essence,
our approach as well. Given the success of the approach in WSMP, MUMPS, and STRUMPACK and the similarity of the 
multifrontal or nested dissection algorithm to the build stage of the HPS algorithm, we should
expect the concept to be successful in accelerating the HPS solver.

\subsection{Outline}
This paper begins with a description of the HPS method in Section \ref{sec:HPS}.  Techniques for finding 
the optimal shared-memory parallelization are presented in Section \ref{sec:optimization}.
Section \ref{sec:results} illustrates the results of the optimization procedure
and the speedup obtained when the HPS method is implemented on a desktop computer.
Finally, the paper closes with remarks and future directions
in Section \ref{sec:summary}.

\section{The HPS method}
\label{sec:HPS}
This section reviews the HPS method for solving the boundary value problem (\ref{eq:basic}).  The solution technique begins by partitioning the domain $\Omega$ into a 
collection of square (or possibly rectangular) boxes, called \textit{leaf boxes}.
Throughout this paper, we denote the parameter for the order of the 
discretization, corresponding to the number of points in each direction on a leaf, as $n_c$.  
For a uniform discretization, the size of all leaf boxes is chosen so that 
any potential solution $u$ of equation (\ref{eq:basic}), as well as its first 
and second derivatives, can be accurately interpolated from their values 
at the local discretization points on any leaf box.

\begin{figure}
\setlength{\unitlength}{1mm}
\begin{picture}(120,40)
\put(05,00){
\includegraphics[width=120mm]{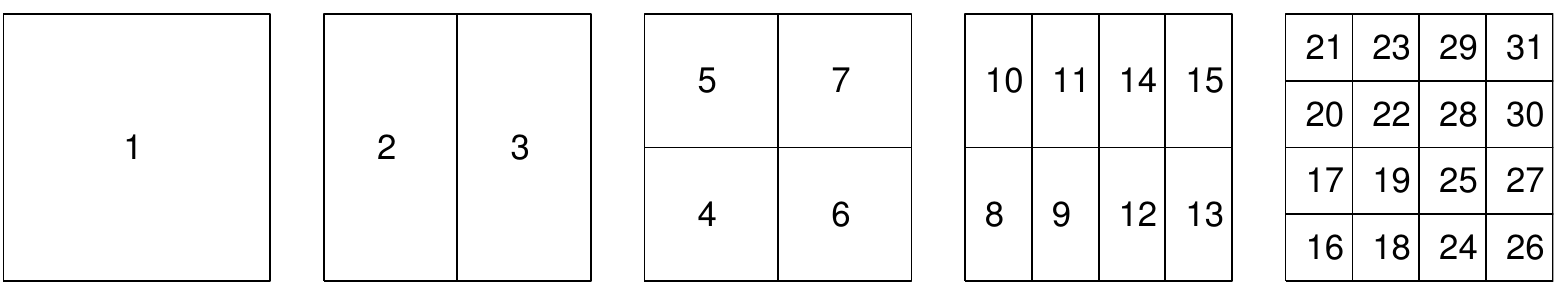}}
\end{picture}
\caption{\label{fig:tree_numbering}
The square domain $\Omega$ is split into $4 \times 4$ leaf boxes.
These are then gathered into a binary tree of successively larger boxes
as described in Section \ref{sec:HPS}. One possible enumeration
of the boxes in the tree is shown, but note that the only restriction is
that if box $\tau$ is the parent of box $\sigma$, then $\tau < \sigma$.}
\end{figure}

Next, a binary tree on the collection of leaf boxes is constructed by
hierarchically merging them. All boxes on
the same level of the tree are roughly of the same size, as shown in Figure
\ref{fig:tree_numbering}.  The boxes should be ordered so
that if $\tau$ is a parent of a box $\sigma$, then $\tau < \sigma$. We
also assume that the root of the tree (i.e.~the full box $\Omega$) has
index $\tau=1$. We let $\Omega^{\tau}$ denote the domain associated with box $\tau$.
Let $N_{\rm boxes}$ denote the number of boxes in the tree.  For the 
tree in Figure \ref{fig:tree_numbering}, $N_{\rm boxes} = 31$.

Recall, from Section 
\ref{sec:overview}, that the solution technique is comprised of 
a \textit{precomputation} stage and a \textit{solve} stage.  
The precomputation stage discretizes the partial differential equation
and builds a direct solver.  The solve stage uses the precomputed
direct solver information applied to body load $s(\vct{x})$
and boundary condition $t(\vct{x})$ information to construct an 
approximate solution the partial differential equation.  
The two major computational components of these stages involve
\textit{leaf} and \textit{merging} computations.

The key to merging boxes is a Poincar\'e-Steklov operator.  
For variable coefficient Helmholtz
problems such as (\ref{eq:basic}), the impedance-to-impedance (ItI) operator is used.  The ItI operator is defined as follows:

\begin{definition}[impedance-to-impedance map]
Fix $\eta\in\mathbb{C}$, and $\Re(\eta) \neq 0$. Let
\begin{eqnarray}
f&:=& u_n+i\eta u|_\Gamma
\label{f}
\\
g&:=&u_n-i\eta u|_\Gamma
\label{g}
\end{eqnarray}
be Robin traces of $u$. We refer to $f$ and $g$ as the ``incoming'' and ``outgoing'' (respectively)
impedance data.
For any $\kappa>0$, the \textit{ItI operator} $R:L^2(\Gamma)\to L^2(\Gamma)$ is defined by
\begin{equation}
R f = g,
\label{R}
\end{equation}
for $f$ and $g$ the Robin traces of $u$ the solution of (\ref{eq:basic}),
for all $f \in L^2(\Gamma)$.
 
\end{definition}

\begin{remark}
For the impedance boundary value problem, the parameter $\eta$ in the definition of 
the ItI operator is the same as the $\eta$ in equation (\ref{eq:basic}).  
 For Dirichlet and Neumann boundary value problems, $\eta = \kappa$ is 
 typically chosen in practice.
\end{remark}

When $s(\vct{x})$ in (\ref{eq:basic}) is nonzero, it is advantageous to 
write the solution $u(\vct{x})$ as the superposition of the 
homogeneous solution $\hmgsoln$
and the particular solution $\prtsoln$; i.e. $u = \hmgsoln+\prtsoln$ where
$\prtsoln$ is the solution of 
\begin{align}
    -\Delta \prtsoln - \kappa^2c(\pxx)\prtsoln &= s(\vct{x}) \qquad & \pxx\in\Omega      \label{eq:prtsolneqn}\\
 \frac{\partial \prtsoln}{\partial \vct{\nu}} +\textrm{i}\eta \prtsoln &= 0 \qquad & \pxx\in \Gamma      \nonumber
\end{align}
and $\hmgsoln$ is the solution of 
\begin{align}
    -\Delta \hmgsoln - \kappa^2c(\pxx)\hmgsoln &= 0 \qquad & \pxx\in\Omega     \label{eq:hmgsolneqn}\\ 
 \frac{\partial \hmgsoln}{\partial \vct{\nu}} +\textit{i}\eta \hmgsoln &= t(\vct{x}) \qquad & \pxx\in  \Gamma.      \nonumber
\end{align}

The remainder of the section presents the technique for 
discretizing leaf boxes (Section \ref{sec:leaf}) and merging neighboring boxes (Section \ref{sec:merge})
via this superposition form.  Specifically, a collection
of approximate solution, impedance, and ItI operators $\mtx{R}$ 
are constructed for each box.  
% presents a technique for discretizing (\ref{eq:basic})
%and constructing an approximate ItI operator $\mtx{R}$ via a 
%modified spectral collocation method.  Then section \ref{sec:merge} presents 
%the technique for merging two boxes constructing corresponding operators. 
An approximate solution can then be constructed for any body load $s(\vct{x})$ 
and boundary condition $t(\vct{x})$ by sweeping the tree twice.  First,
particular solution information is constructed, moving from the
leaf boxes up the hierarchical tree.  The approximate solution is then
created by propagating homogeneous boundary information down the 
tree.  Algorithm 2 presents the details of this procedure.  When there 
is no body load (i.e. $s(\vct{x}) = 0$), the solution procedure needs only
 the downward sweep of the tree.  The homogeneous solver is the 
same algorithm as presented in \cite{2013_martinsson_ItI}.  

\subsection{Leaf computation}
\label{sec:leaf}
This section presents a modified spectral collocation method for
constructing the necessary operators for processing a leaf box $\tau$.  
The modified spectral collocation 
technique, first presented in \cite{2018_geldermans_gillman_adaptive}, is 
ideal for the HPS method because it does not involve corner discretization 
points, for which Poincar\'e-Steklov operators are not always well defined.

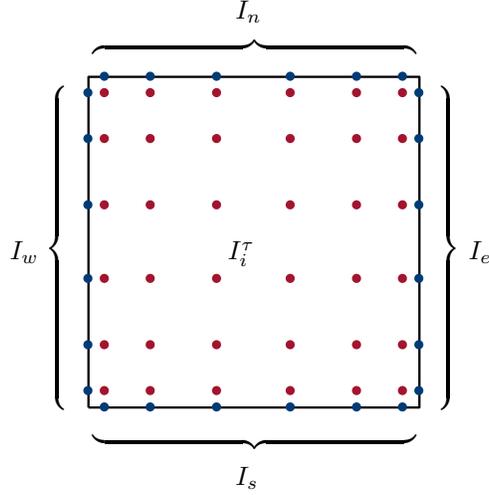
\begin{figure}  
\centering
\setlength{\unitlength}{1mm}
\begin{picture}(60,60)
\put(05,05){\includegraphics[height = 50mm]{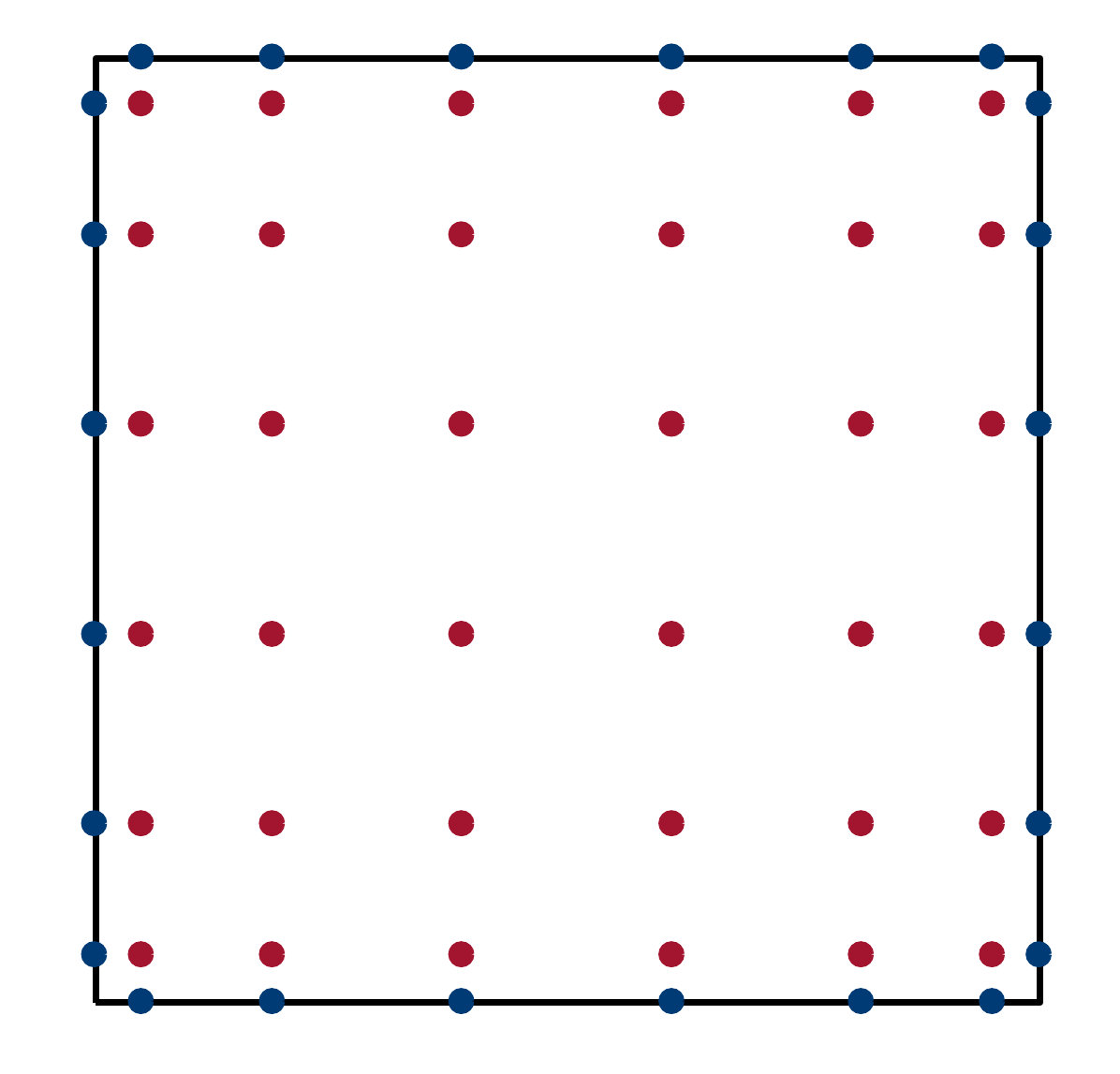}}
\put(28,28){$I^\tau_i$}
\put(10,05){$\underbrace{\hspace{4.3cm}}$}
\put(29,-2){$I_s$}
\put(10,55){$\overbrace{\hspace{4.3cm}}$}
\put(29,60){$I_n$}
\put(56,08){{\rotatebox{90}{$\underbrace{\hspace{4.3cm}}$}}}
\put(60,28){$I_e$}
\put(4,08){{\rotatebox{90}{$\overbrace{\hspace{4.3cm}}$}}}
\put(-1,28){$I_w$}
\end{picture}
\caption{\label{fig:leaf} Illustration of the discretization 
points for a leaf box $\tau$.  The points in blue are the 
boundary points with indices $I^\tau_b = [I_s,I_e,I_n,I_w]$.
The points in red are the interior points with indices $I^\tau_i$.  
The points in black are the omitted corner points.
}
\end{figure}

The modified spectral collocation technique begins with the 
classic $n_c\times n_c$ product Chebychev grid and the 
corresponding standard spectral differential matrices $\mtx{D}_x$ and $\mtx{D}_y$, as defined in 
 \cite{2000_trefethen_spectral_matlab}.   
Let ${I}^\tau_i$ denote the index vector corresponding to points on the 
interior of $\Omega^\tau$ and ${I}^\tau_b$ denote the index vector 
corresponding to points on the boundary of $\Omega^\tau$, \textbf{not} 
including the corner points, based on the classic tensor grid.
Figure \ref{fig:leaf} illustrates the 
indexing of the points in terms of the classic discretization.  Thus 
$\{\pxx_j\}_{j=1}^{n_c^2-4}$ denotes the discretization points 
in $\Omega^\tau$ given by the union of the red and blue points 
in Figure \ref{fig:leaf}.  
We order the solution vector $\vct{u}$ according to the following:
$\vct{u} = \vtwo{\vct{u}_b}{\vct{u}_i}$ where $\vct{u}_b$ and $\vct{u}_i$ 
denote the approximate values of the solution on the 
boundary and the interior, respectively.  The homogeneous and particular 
solution vectors $\hmgsolnvec$ and $\prtsolnvec$ are ordered in a similar
manner.  
The ordering of the 
entries related to the boundary corresponding to the discretization points 
is $I^\tau_b = [I_s,I_e,I_n,I_w]$ where $I_s$ denotes the blue points on the south
boundary in Figure \ref{fig:leaf}, etc. Let $I^\tau = [I^\tau_b,I^\tau_i]$ denote
the collection of all indices that are used in the discretization.  Let 
the $n_b= 4n_c-8$ denote the length of $I^\tau_b$ and $n_i=(n_c-2)^2$ denote the length of $I^\tau_i$.
The number of discretization points on the leaf box $\tau$ is $n^\tau = n_b+n_i$.

The discrete approximation of the differential operator on $\Omega^\tau$ using 
classic spectral collocation \cite{2000_trefethen_spectral_matlab} is given by 
% $$\mtx{A} = -\mtx{D}_x^2-\mtx{D}^2_y+{\rm{diag}}\left\{b(\pxx_j)\right\}.$$
\begin{equation*}
\mtx{\hat{A}} =
-\mtx{D}_x^{2}-\mtx{D}_y^{2}
-\mtx{C},
\end{equation*}
where $\mtx{C}$ is the diagonal matrix with diagonal entries $\{\kappa^2c(\pxx_{k})\}_{k=1}^{n_c^{2}}$.

Due to the tensor product basis, we know the entries of $\mtx{D}_x$ and 
$\mtx{D}_y$ corresponding to the interaction of the corner points with 
the points on the interior of $\Omega^\tau$ are zero.  The directional 
basis functions for the other points on the boundary are not impacted by 
the removal of the corner points.  Thus the differential operators from
the classic pseudospectral discretization can be used to create the
approximation of the local differential operator, the ItI operator, and 
all other necessary leaf operators.

This information allows for the approximation of the differential operator on
$\Omega^\tau$ using the modified discretization to be constructed from the 
classic spectral collocation differential operator.  Specifically, the 
approximate modified spectral collocation operator is 
the $n^\tau\times n^\tau$ matrix 
$$ \mtx{A} = \mtx{\hat{A}}(I^\tau,I^\tau)$$
where $n^\tau = n_c^2-4$.

Likewise, operators can be constructed to approximate impedance operators. 
Let $\mtx{N}$ denote the $n_b\times n^\tau$ matrix that takes normal derivatives of the basis 
functions.  Then $\mtx{N}$ is given by 

$$\mtx{N} = \left[\begin{array}{r} -\mtx{D}_x(I_s,I^\tau) \\
                           \mtx{D}_y(I_e,I^\tau)\\
                           \mtx{D}_x(I_n,I^\tau)\\
                           -\mtx{D}_y(I_w,I^\tau)\end{array}\right].$$              
The $n_b\times n^\tau$ matrix for creating the incoming impedance data is 
$$\mtx{F} = \mtx{N}+\textrm{i}\eta\mtx{I}_{n_c^2}(I^\tau_b,I^\tau)$$
and the $n_b\times n^\tau$ matrix for creating the outgoing impedance data is 
$$\mtx{G} = \mtx{N}-\textrm{i}\eta\mtx{I}_{n_c^2}(I^\tau_b,I^\tau)$$
where $\mtx{I}_{n_c^2}$ is the identity matrix of size 
$n_c^2$.

\subsubsection{Homogeneous solution operators}
\label{sec:homo_leaf}
To construct the homogeneous solution operators, we consider 
the discretized differential equation defined on $\Omega^\tau$.
The discretized body-load problem, to find the approximation to $\hmgsoln$ at the 
collocation points takes, the form

\begin{equation}
 \mtx{B}\vtwo{\hmgsolnvec_{b}}{\hmgsolnvec_i} = 
 \vtwo{\mtx{F}}{\mtx{\hat{A}}(i,b) \ \mtx{\hat{A}}(i,i)} \vtwo{\hmgsolnvec_{b}}{\hmgsolnvec_i} = \vtwo{\mtx{t}}{\vct{0}},
 \label{eq:homogeneous}
 \end{equation}
where $\mtx{B}$ is an $n^\tau \times n^\tau$ matrix,
$\hmgsolnvec$ is the vector with the approximate homogeneous solution 
at the collocation points, $\mtx{\hat{A}}_{i,i} = \mtx{\hat{A}}({I}^\tau_i,{I}^\tau_i)$ is a matrix of size
${n_i \times n_i}$, 
$\mtx{\hat{A}}_{i,b} = \mtx{\hat{A}}({I}^\tau_i,{I}^\tau_b)$ is a matrix of size $n_i \times n_b$, 
and $\vct{t}$ is vector of length $n_b$ containing impedance boundary data.

The \textit{homogeneous solution operator} $\HmgSolnOp^\tau$ is the $n^\tau \times n_b$ matrix
defined by solving 
\begin{equation}
\mtx{B}\HmgSolnOp^\tau= \vtwo{\vct{I}_{n_b}}{\vct{0}_{ n_i \times n_b}}.
\label{eq:hmgsolnopdef}
\end{equation}

To construct the approximate \textit{ItI operator} $\mtx{R}^\tau$ (of size
$n_b\times n_b$), we simply need 
to apply $\mtx{G}$ to $\HmgSolnOp^\tau$, that is 
$$\mtx{R}^\tau = \mtx{G}\HmgSolnOp^\tau.$$

\subsubsection{Particular solution operators}
\label{sec:part_leaf}
The particular solution operators are constructed in a 
similar manner.  Specifically, the discretized version of 
(\ref{eq:prtsolneqn}) takes the form
\begin{equation}
 \mtx{B}\vtwo{\vct{\prtsoln}_{b}}{\vct{\prtsoln}_i} = 
 \vtwo{\mtx{F}}{\mtx{\hat{A}}(i,i) \ \mtx{\hat{A}}(i,b)} \vct{\prtsoln} = \vtwo{\vct{0}}{\mtx{s}},
 \label{eq:particular}
 \end{equation}
 where $\mtx{B}$ is an $n^\tau \times n^\tau$ matrix, and $\mtx{s}$ is a vector of length $n_i$
containing body load data.  
 
The \textit{particular solution operator} $\PrtSolnOp^\tau$ is an $n^\tau \times n_i$ matrix  
which can be used to approximate the particular solution $\prtsoln$ 
on the leaf $\tau$ when applied to any body load vector $\mtx{s}$. It is the solution of 
\begin{equation}
\mtx{B}\PrtSolnOp^\tau = \vtwo{\vct{0}_{n_b\times n_i}}{\vct{I}_{n_i}}.
\label{eq:prtsolnopdef}
\end{equation}

As with the homogeneous case, the operator constructing the approximation
of the outgoing impedance data is constructed by applying the operator 
$\mtx{G}$ to the particular solution operator $\PrtSolnOp^\tau$.  Let 
the $n_b\times n_i$ matrix 
$$\mtx{\Gamma}^\tau = \mtx{G}\PrtSolnOp^\tau$$
denote this \textit{particular solution outgoing impedance} operator.

\begin{remark}
 Once all the leaf operators are constructed for a box $\tau$,
 the solution vector $\vct{u}^\tau$ is given by 
 $$\vct{u}^\tau = \HmgSolnOp^\tau \mtx{t}+\PrtSolnOp^\tau \vct{s} = 
 \hmgsolnvec^\tau + \prtsolnvec^\tau,$$
 where $\mtx{t}$ is a vector whose entries are impedance boundary 
 data at the boundary nodes on $\tau$ and $\vct{s}$ is a 
 vector whose entries are the evaluation of the body load $s(\vct{x})$ at the interior discretization points of $\tau$.  
 The outgoing impedance data is given by 
 $$\vct{g}^\tau = \mtx{R}^\tau \mtx{t} + \mtx{\Gamma}^\tau \vct{s}
  = \mtx{R}^\tau \mtx{t} + \vct{h}^\tau,$$
  where $\vct{h}^\tau$ is the \textit{particular solution outgoing impedance} boundary data. 
\end{remark}

%  We first turn to the construction of the homogeneous and particular solution operators
%  on the leaf boxes, as well as the ItI operator.  The computation uses a modified spectral
%  collocation method first described in \cite{2018_geldermans_gillman_adaptive}.  The key
%  modification is the removal of the corner points from the classical two-dimensional 
%  tensor product basis formed from the Chebyshev points, illustrated in Figure \ref{fig:leaf}.  Removal of the corner points
%  improves the stability of the method, with minimal impact on the accuracy because
%  the corner points do not belong to the stencil of the discretized differential operator for the interior points.
% 

\subsection{Merging two boxes}
\label{sec:merge}
This section presents the technique for constructing the necessary 
operators for the union of two boxes for which outgoing particular
solution information and ItI operators have already been constructed.

\begin{figure}
\centering
\setlength{\unitlength}{1mm}
\begin{picture}(95,55)
 \put(00,00){\includegraphics[trim=1.1in 3.8in 0.9in 3.5in, clip, height=55mm]{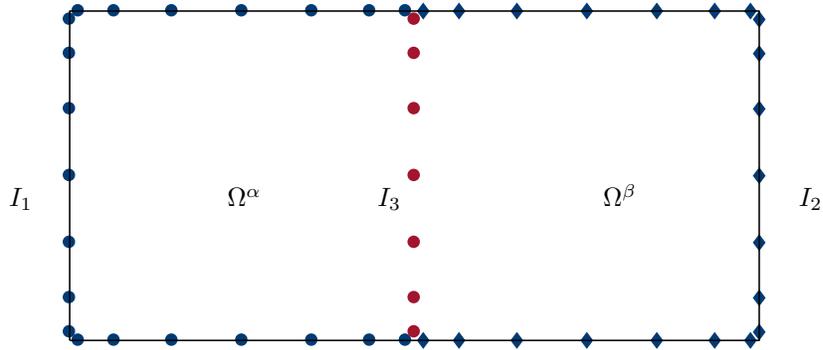}}
\put(24,25){$\Omega^{\alpha}$}
\put(74,25){$\Omega^{\beta}$}
\put(-5,25){$I_{1}$}
\put(100,25){$I_{2}$}
\put(44,25){$I_{3}$}
% \put(47,52){$I_{3}$}
\end{picture}
\caption{Notation for the merge operation described in Section \ref{sec:merge}.
The rectangular domain $\Omega$ is formed by two squares $\Omega^{\alpha}$ 
and $\Omega^{\beta}$.  The sets $I_{1}$ (blue circles) and $I_{2}$ (blue diamonds) 
form the exterior nodes, while $I_{3}$ (red circles) consists of the interior nodes.}
\label{fig:siblings_notation}
\end{figure}

Let $\Omega^{\tau}$ denote a box with children $\Omega^{\alpha}$ and
$\Omega^{\beta}$ so that
$$
\Omega^{\tau} = \Omega^{\alpha} \cup \Omega^{\beta}.$$ 
For concreteness, but without loss of generality, assume that $\Omega^{\alpha}$ and
$\Omega^{\beta}$ share a vertical edge as shown in Figure \ref{fig:siblings_notation}.  
We partition the points on $\partial\Omega^{\alpha}$ and $\partial\Omega^{\beta}$ into three sets:
\begin{tabbing}
\mbox{}\hspace{5mm}\= $I_{1}$ \hspace{4mm} \=
Boundary nodes of $\Omega^{\alpha}$ that are not boundary nodes of $\Omega^{\beta}$.\\
\> $I_{2}$ \> Boundary nodes of $\Omega^{\beta}$ that are not boundary nodes of $\Omega^{\alpha}$.\\
\> $I_{3}$ \> Boundary nodes of both $\Omega^{\alpha}$ and $\Omega^{\beta}$ that are \textit{not} boundary nodes of the\\
\> \>  union box $\Omega^{\tau}$.
\end{tabbing}
The indexing for the points on the interior and boundary of 
$\Omega^\tau$ are $I_{\rm i}^\tau = I_3$ and $I_{\rm b}^\tau = [I_1,I_2]$, 
respectively.

For the box $\alpha$, let $\vct{t}^\alpha$ denote the \textit{homogeneous solution} incoming
impedance boundary data, $\vct{h}^\alpha$ denote the \textit{particular solution} outgoing
impedance boundary data, and $\vct{g}^\alpha$ denote the \textit{total} outgoing impedance
boundary data.  Define the vectors $\vct{t}^\beta$, $\vct{h}^\beta$, and $\vct{g}^\beta$
similarly.  

Using the ItI operators $\mtx{R}^\alpha$ and $\mtx{R}^\beta$ and ordering 
everything according to the boundary numbering in Figure \ref{fig:siblings_notation}, the 
outgoing impedance data for boxes $\alpha$ and $\beta$ are given by 

\begin{equation}  \vtwo{\vct{g}^\alpha_1}{\vct{g}^\alpha_3}=\mtwo{\mtx{R}^\alpha_{11}}{\mtx{R}^\alpha_{13}}{\mtx{R}^\alpha_{31}}{\mtx{R}^\alpha_{33}} 
\vtwo{\vct{t}^\alpha_1}{\vct{t}^\alpha_3} 
+ \vtwo{\vct{h}^\alpha_1}{\vct{h}^\alpha_3} 
\label{eq:ItIalp}\end{equation}
and 
\begin{equation}\vtwo{\vct{g}^\beta_2}{\vct{g}^\beta_3} = 
\mtwo{\mtx{R}^\beta_{22}}{\mtx{R}^\beta_{23}}{\mtx{R}^\beta_{32}}{\mtx{R}^\beta_{33}} 
\vtwo{\vct{t}^\beta_2}{\vct{t}^\beta_3}
 + \vtwo{\vct{h}^\beta_2}{\vct{h}^\beta_3}.
 \label{eq:ItIbet}\end{equation}

% where $\vtwo{\vct{h}^\alpha_1}{\vct{h}^\alpha_3}$ and $\vtwo{\vct{h}^\beta_2}{\vct{h}^\beta_3}$
% are the outgoing impedance data due to the particular solutions on each box.

Since the normal vectors are opposite in each box, we know $\vct{t}^\alpha_3 = -\vct{g}^\beta_3$ and $\vct{g}^\alpha_3 = -\vct{t}^\beta_3$.
Using this information in the bottom row equations in (\ref{eq:ItIalp})
and (\ref{eq:ItIbet}), $\vct{t}^\alpha_3$ and $\vct{t}^\beta_3$ can 
found in terms of $\vct{t}_1^\alpha$, $\vct{t}_2^\beta$, $\vct{h}^3_\alpha$, and $\vct{h}^3_\beta$.
They are given by
\begin{equation}
\vct{t}_3^\alpha = \mtx{\Phi}^\alpha \vtwo{\vct{t}_1^\alpha}{\vct{t}_{2}^\beta} + \mtx{\Upsilon}^\alpha \vtwo{\vct{h}^\alpha_3}{\vct{h}_3^\beta}
%\begin{equation}
%\vct{t}_3^\alpha = \underbrace{\mtx{W}^{-1}\left[\mtx{R}_{33}^\beta \mtx{R}_{31}^\alpha | -\mtx{R}_{32}^\beta\right] }_{\mtx{\Phi}^\alpha}
%\vtwo{\vct{t}_1^\alpha}{\vct{t}_{2}^\beta} +\underbrace{\left[
%\mtx{W}^{-1}\mtx{R}_{33}^\beta\ | \ -\mtx{W}^{-1}\right]}_{\mtx{\Upsilon}^\alpha}\vtwo{\vct{h}^\alpha_3}{\vct{h}_3^\beta}
 \label{eq:talpha}
\end{equation}
and 

\begin{equation}
\vct{t}^\beta_3 = {\mtx{\Phi}^\beta}
\vtwo{\vct{t}_{1}^\alpha}{\vct{t}_{2}^\beta} + {\mtx{\Upsilon}^\beta}\vtwo{\vct{h}^\alpha_3}{\vct{h}_3^\beta},
 \label{eq:tbeta}
\end{equation}

where 
\begin{align*}
\mtx{\Phi}^\alpha &= \mtx{W}^{-1}\left[\mtx{R}_{33}^\beta \mtx{R}_{31}^\alpha | -\mtx{R}_{32}^\beta\right],\\
\mtx{\Upsilon}^\alpha &= \left[\mtx{W}^{-1}\mtx{R}_{33}^\beta\ | \ -\mtx{W}^{-1}\right],\\
\mtx{\Phi}^\beta &= \left[-\mtx{R}_{31}^\alpha -\mtx{R}^\alpha_{33}\mtx{W}^{-1}\mtx{R}_{33}^\beta\mtx{R}^\alpha_{31} | \mtx{R}_{33}^\alpha \mtx{W}^{-1}\mtx{R}_{32}^\beta\right],\\
 \mtx{\Upsilon}^\beta &= \left[-\left(\mtx{I}+\mtx{R}^\alpha_{33}\mtx{W}^{-1}\mtx{R}_{33}^\beta\right)  \ | \ \mtx{R}_{33}^\alpha \mtx{W}^{-1}\right], \ {\rm and}\\
\mtx{W} &= \mtx{I}-\mtx{R}_{33}^\beta\mtx{R}^\alpha_{33}.
 \end{align*}

Substituting (\ref{eq:talpha}) and (\ref{eq:tbeta}) into the top row equations of (\ref{eq:ItIalp}) and (\ref{eq:ItIbet}) results in the following expression 
for the outgoing impedance data for the box $\Omega^\alpha\cup \Omega^\beta$

\begin{equation}
 \vtwo{\vct{g}^\alpha_1}{\vct{g}^\beta_2} = \mtx{R}^\tau
\vtwo{\vct{t}_1^\alpha}{\vct{t}_2^\beta}  +\vtwo{\vct{h}^\alpha_1}{\vct{h}_2^\beta} +
\mtx{\Gamma}^\tau
 \vtwo{\vct{h}^\alpha_3}{\vct{h}^\beta_3},
 \label{eq:out}
\end{equation}
where 
$$\mtx{R}^\tau= \mtwo{\mtx{R}_{11}^\alpha}{\mtx{0}}{\mtx{0}}{\mtx{R}_{22}^\beta} + \mtwo{\mtx{R}^\alpha_{13}}{\mtx{0}}{\mtx{0}}{\mtx{R}^\beta_{23}} \vtwo{\mtx{\Phi}^\alpha}{\mtx{\Phi}^\beta}$$ 
is the homogeneous ItI operator
and 
$$\mtx{\Gamma}^\tau = \mtwo{\mtx{R}^\alpha_{13}}{\mtx{0}}{\mtx{0}}{\mtx{R}^\beta_{23}} \vtwo{\mtx{\Upsilon}^\alpha}{\mtx{\Upsilon}^\beta} $$ is the outgoing particular solution flux due to interior edge operator.  The outgoing
particular solution flux is  
\begin{equation}
\vct{h}^\tau = \vtwo{\vct{h}^\alpha_1}{\vct{h}_2^\beta} +
\mtx{\Gamma}^\tau
 \vtwo{\vct{h}^\alpha_3}{\vct{h}^\beta_3}.
 \label{eq:flux_part}
 \end{equation}

\begin{remark}
When $\tau = 1$, the boundary data from equation (\ref{eq:basic}) gets utilized 
in equation (\ref{eq:out}) at the $\vct{t}$ contributions.  The $\vct{h}$ 
contributions is determined from $\alpha = 2$ and $\beta = 3$ boxes.   
\end{remark}
 
% % %  \begin{remark}
% % %  \label{re:store}
% % %  For computational efficiency, the operators $\mtx{\Upsilon}^\alpha$, $\mtx{\Upsilon}^\beta$, 
% % %  and $\mtx{\Gamma}^\tau$ are not explicitly formed.
% % %  Instead, the matrices $\mtx{W}^{-1}$, $\mtx{R}_{33}^\alpha$, $\mtx{R}_{33}^\beta$, $\mtx{R}_{13}^\alpha$, and $\mtx{R}_{23}^\beta$
% % %  needed to catpure of the action of the operators are stored.
% % %  \end{remark}

\subsection{Computational cost}
\label{sec:cost}
The cost of constructing the discretization and direct
solver is dominated by inverting the matrix $\mtx{W}$ 
of size $O(N^{1/2})$ at the top level in the tree, where $N$ is the number of discretization points.  Thus
the computational cost of the precomputation stage is 
$O(N^{3/2})$.
At any level in the solve stage, the cost of applying all the 
operators is $O(N)$ and there are $\log N$ levels in a uniform 
tree.  It follows that the total cost of the apply the solver is $O(N\log N)$ with a small constant.

\subsection{Operator storage}
\label{sec:mem}
Most operators are explicitly computed during the build stage and stored.  This allows
for the solve stage to simply be a collection of matrix vector multiplies.  
Specifically, for each leaf box, the homogeneous
and particular solution operators $\mtx{\Psi}$ and $\mtx{Y}$, as well as the operator $\mtx{\Gamma}$ for the outgoing
impedance associated with the particular solution are stored.  
For boxes that are processed via the merge procedure, the build stage is more efficient if 
all the operators are not stored directly. The operators $\mtx{\Phi}^\alpha$ 
and $\mtx{\Phi}^\beta$ for computing the homogeneous solution incoming impedance at the interface are computed and stored.  However, rather than compute $\mtx{\Upsilon}^\alpha$, $\mtx{\Upsilon}^\beta$, and $\mtx{\Gamma}^\tau$
for each parent box, the set of matrices that capture the action of applying these operators (namely $\mtx{W}^{-1}$, $\mtx{R}_{33}^\alpha$, $\mtx{R}_{33}^\beta$, $\mtx{R}_{13}^\alpha$,
and $\mtx{R}_{23}^\beta$) are stored.
The matrices $\mtx{W}^{-1}$, $\mtx{R}_{33}^\alpha$, $\mtx{R}_{33}^\beta$, $\mtx{R}_{13}^\alpha$,
and $\mtx{R}_{23}^\beta$ create the action 
of  $\mtx{\Upsilon}^\alpha$, $\mtx{\Upsilon}^\beta$, and $\mtx{\Gamma}^\tau$ in the upward pass upward pass of the solve stage via only a series of matrix-vector multiplications.  This allows for the computation of
the particular solution information to be computed for 
less (total) computational cost since matrix-vector products are more computationally efficient than the 
 matrix-matrix products required to explicitly build the $\mtx{\Upsilon}$ and $\mtx{\Gamma}$ operators.  This efficiency comes at a 
cost of storing an additional matrix of size $n_3^\tau \times n_3^\tau$ where $n_3^\tau$ is the length
of $I_3$.  The ItI operator $\mtx{R}$ of all children boxes is deleted once the parent box is processed.

\begin{figure}[h!]
\fbox{
\begin{minipage}{115mm}
\begin{center}
\textsc{Algorithm 1} (Precomputation stage)
\end{center}

This algorithm builds all the operator needed to construct
an approximate solution to (\ref{eq:basic}) for any choice of body load $s(\vct{x})$
and incoming impedance boundary condition $t(\vct{x})$.
It is assumed that if node $\tau$ is a parent of node $\sigma$, then $\tau < \sigma$.
Let $N_{\rm boxes}$ denote the number of boxes in the tree.

\rule{\textwidth}{0.5pt}

\begin{tabbing}
\mbox{}\hspace{7mm} \= \mbox{}\hspace{6mm} \= \mbox{}\hspace{6mm} \= \mbox{}\hspace{6mm} \= \mbox{}\hspace{6mm} \= \kill
(1)\> \textbf{for} $\tau = N_{\rm boxes},\,N_{\rm boxes}-1,\,N_{\rm boxes}-2,\,\dots,\,1$\\
(2)\> \> \textbf{if} ($\tau$ is a leaf)\\
(3)\> \> \> Construct $\mtx{R}^{\tau}$, $\mtx{Y}^\tau$, $\mtx{\Psi}^\tau$ and $\mtx{\Gamma}^\tau$ via the process described in\\
\> \> \> Section \ref{sec:leaf}.\\
(4)\> \> \textbf{else}\\
(5)\> \> \> Let $\alpha$ and $\beta$ be the children of $\tau$.\\
(6)\> \> \> Split $I_{\rm b}^{\alpha}$ and $I_{\rm b}^{\beta}$ into vectors $I_{1}$, $I_{2}$, and $I_{3}$ as shown in \\
\> \> \> Figure \ref{fig:siblings_notation}.\\
(7)\> \> \> $\mtx{W} = \mtx{I}-\mtx{R}^{\beta}_{33}\mtx{R}^{\alpha}_{33}$\\
(8)\> \> \> $\mtx{\Phi}^{\alpha} = \mtx{W}^{-1}\left[\mtx{R}_{33}^{\beta} \mtx{R}_{31}^{\alpha} | -\mtx{R}_{32}^{\beta}\right]$\\
(9) \> \> \> $\mtx{\Upsilon}^{\alpha} = \left[\mtx{W}^{-1}\mtx{R}_{33}^{\beta}\ | \ -\mtx{W}^{-1}\right] $\\
(10)\> \> \> $\mtx{\Phi}^{\beta} = \left[-\mtx{R}_{31}^{\alpha} -\mtx{R}^{\alpha}_{33}\mtx{W}^{-1}\mtx{R}_{33}^{\beta}\mtx{R}^{\alpha}_{31} | \mtx{R}_{33}^{\alpha} \mtx{W}^{-1}\mtx{R}_{32}^{\beta}\right] $\\ 
(11)\> \> \> $\mtx{\Upsilon}^{\beta} = \left[-\left(\mtx{I}+\mtx{R}^{\alpha}_{33}\mtx{W}^{-1}\mtx{R}_{33}^{\beta}\right)  \ | \ \mtx{R}_{33}^{\alpha} \mtx{W}^{-1}\right] $\\  
(12)\> \> \> $\mtx{R}^\tau = \mtwo{\mtx{R}_{11}^{\alpha}}{\mtx{0}}{\mtx{0}}{\mtx{R}_{22}^{\beta}} + \vtwo{\mtx{R}^{\alpha}_{13}}{\mtx{R}^{\beta}_{23}} \vtwo{\mtx{\Phi}^{\alpha}}{\mtx{\Phi}^{\beta}}$\\
(13)\> \> \> $\mtx{\Gamma}^\tau = \vtwo{\mtx{R}^{\alpha}_{13}}{\mtx{R}^{\beta}_{23}} \vtwo{\mtx{\Upsilon}^{\alpha}}{\mtx{\Upsilon}^{\beta}}$\\
(14)\> \> \> Delete $\mtx{R}^{\alpha}$ and $\mtx{R}^{\beta}$.\\
(15)\> \> \textbf{end if}\\
(16)\> \textbf{end for}
\end{tabbing}
\end{minipage}}
\label{alg:precompute}
\end{figure}

\begin{figure}[h!]
\fbox{
\begin{minipage}{115mm}
\begin{center}
\textsc{Algorithm 2} (Solve stage)
\end{center}
This algorithm constructs an approximate solution $\uu$
to (\ref{eq:basic}) given the body load $s(\vct{x})$ 
and incoming impedance boundary condition $t(\vct{x})$.
It is assumed that if node $\tau$ is a parent of node $\sigma$, then $\tau < \sigma$.
Let $N_{\rm boxes}$ denote the number of boxes in the tree.
All operators are assumed to be precomputed.

\rule{\textwidth}{0.5pt}
\begin{tabbing}
\mbox{}\hspace{7mm} \= \mbox{}\hspace{6mm} \= \mbox{}\hspace{6mm} \= \mbox{}\hspace{6mm} \= \mbox{}\hspace{6mm} \= \kill
Upward pass\\
(1)\> \textbf{for} $\tau = N_{\rm boxes},\,N_{\rm boxes}-1,\,N_{\rm boxes}-2,\,\dots,\,1$\\
(2)\> \> \textbf{if} ($\tau$ is a leaf)\\
(3)\> \> \> Compute $\prtsolnvec^\tau = \PrtSolnOp^\tau\vct{s}$ and $\outimpprtvec^\tau = \OutImpPrtSolnOp^\tau\vct{s}$ for the leaf.\\
(4)\> \> \textbf{else}\\
(5)\> \> \> Let $\alpha$ and $\beta$ be the children of $\tau$.\\
(6)\> \> \> Compute $\tilde{\vct{t}}^{\alpha}_3 = \mtx{\Upsilon}^{\alpha}\vtwo{\vct{h}^{\alpha}_3}{\vct{h}^{\beta}_3}$ and $\tilde{\vct{t}}^{\beta}_3 = \mtx{\Upsilon}^{\beta}\vtwo{\vct{h}^{\alpha}_3}{\vct{h}^{\beta}_3}$.\\
(7)\> \> \> Compute $\outimpprtvec^\tau$ via (\ref{eq:flux_part}).\\
(8)\> \> \textbf{end if}\\
(9)\> \textbf{end for}\\
Downward pass\\
(10)\> \textbf{for} $\tau = 1,\,2,\,3,\,\dots,\,N_{\rm boxes}$\\
(11)\> \> \textbf{if} ($\tau$ is a leaf)\\
(12)\> \> \> Let $J^\tau$ denote the indices for the 
discretization points in $\tau$.  \\
(13)\> \> \>$\uu(J^\tau) =\HmgSolnOp^{\tau}\,\vct{t}^\tau + \tilde{\uu}^\tau$.\\
(14)\> \> \textbf{else} \\
(15)\> \> \> Let $\alpha$ and $\beta$ be the children of $\tau$.\\
(16)\> \> \>$\vct{t}^{\alpha}_3 = \IncImpHmgSolnOp^{\alpha}\vct{t}^\tau + \tilde{\vct{t}}^{\alpha}$, $\vct{t}^{\beta}_3 = \IncImpHmgSolnOp^{\beta}\vct{t}^\tau + \tilde{\vct{t}}^{\beta}$.\\
(17)\>\> \textbf{end if}\\
(18)\> \textbf{end for}
\end{tabbing}
\end{minipage}}
\label{alg:solve}
\end{figure}

\section{General thread optimization technique}
\label{sec:optimization}
This section presents the proposed optimization technique for shared-memory 
parallelism via OpenMP and the Intel MKL library. Recall from Section \ref{sec:HPS} that the bulk of the computational cost
in both stages of the algorithm is associated with matrix inversion and matrix-matrix multiplication.  

The linear algebra needed to process a given box is essentially sequential. 
For example, in order to construct any of the operators
in the merge, the inverse of $\mtx{W}$ must be formed.  We construct this inverse using MKL, with the routine \verb|zgetrf| computing an LU factorization of $\mtx{W}$ and 
the triangular solve routine  \verb|zgetri| computing the inverse. The dominant computational costs
of the rest of the merge process come from dense matrix multiplications, implemented with 
\verb|zgemm|.  
The operators $\mtx{\Phi}^{\alpha}$ and $\mtx{\Phi}^{\beta}$ are needed for the downward pass of the solve stage and \ also form part of the full ItI operator $\mtx{R}^\tau$ (see section \ref{sec:merge}).  
(As stated in section \ref{sec:mem}, we do not form the particular solution operators explicitly for greater computational efficiency.) 
 Since these operations build on each other and cannot 
be computed simultaneously, exploiting the parallel 
linear algebra capabilities of MKL is 
 important for achieving parallelism within the merge procedure, especially at the top of the tree where the matrices are large. 

The algorithm is domain decomposing, and all boxes on a given level 
in the tree are independent of each other.  As a result, there are two types of 
parallelism that can be exploited: dividing boxes among threads and utilizing the 
multi-threaded linear algebra in MKL.  We propose a hybrid of these the approaches. 
In the bottom of the tree where the matrices involved are small, it is best 
to use a ``divide-and-conquer'' approach which distributes boxes among all 
available threads.  At the top of the tree, it is most efficient to exploit black-box parallel linear algebra routines.  
The best distribution of work on the intermediate levels
depends on the available computational resources and number of boxes on a given level.  This section 
presents a technique for distributing work in the hybrid parallelism setting.

Let $\theta_t$ denote the total number of available threads.  For a level $l$ in the tree, let  $\theta^l_o$ 
denote the number of \textit{outer threads} dedicated to a divide-and-conquer
distribution of boxes, and $\theta^l_i$ denote the number of \textit{inner threads}
given to each outer thread for parallel linear algebra, on that level. Then  $\theta_t \geq \theta^l_o \times \theta^l_i$ and $  \theta^l_o, \theta^l_i \in \mathbb{Z}^+$
for all $l$.  
The implementation of outer and inner threads works as follows.  For each (non-leaf) level $l$ of the tree with $N_\text{boxes}^l$ boxes,  
$2N^l_\text{boxes}$ children boxes are merged by looping over the $N_\text{boxes}^l$  parent boxes.
It is natural to parallelize this loop by placing an OpenMP parallel region with $\theta_o^l$ threads inside of it.  
For a uniform grid, static scheduling is sufficient since the work required to merge each box is identical.  Therefore, the 
best strategy should be to split the boxes as evenly as possible among the threads, which static scheduling does. Within 
this parallel loop, we set the number of threads
available to MKL at $\theta_i^l$ using \verb|mkl_set_num_threads_local|.  Similar adjustments are made for 
the leaf processing as well as the upward and downward pass portions of the solve stage, with their respective outer and inner thread pairs 
for each level. Algorithms 3 and 4 detail the parallel versions of Algorithms 1 and 2.

\begin{remark}If the grid was non-uniform, dynamic scheduling should be used, as the benefit from runtime load balancing would likely outweigh the additional overhead.
\end{remark}

\begin{figure}[h!]
\fbox{
\begin{minipage}{115mm}
\begin{center}
\textsc{Algorithm 3} (Precomputation stage, parallel)
\end{center}

This algorithm presents the parallel version of the 
precomputation stage. The serial version is 
presented in Algorithm 1.

Let $N_\text{boxes}^l$ be 
the number of boxes on a level $l$, and $\theta_o^l$ and $\theta_i^l$ be the outer and inner build stage threads
for a level, respectively.

\rule{\textwidth}{0.5pt}

\begin{tabbing}
\mbox{}\hspace{7mm} \= \mbox{}\hspace{6mm} \= \mbox{}\hspace{6mm} \= \mbox{}\hspace{6mm} \= \mbox{}\hspace{6mm} \= \kill
Leaf construction\\
(1) \> Set outer threads to  $\theta_o^L$. \\
(2) \> \textbf{parfor} $\tau = N_\text{boxes} - N_\text{boxes}^L + 1,  \, N_\text{boxes} - N_\text{boxes}^L + 2, \dots, \, N_\text{boxes}$ \\
(3) \>\> Set inner threads to $\theta_i^L$.\\
(4) \> \> Construct $\mtx{R}^{\tau}$, $\mtx{Y}^\tau$, $\mtx{\Psi}^\tau$ and $\mtx{\Gamma}^\tau$ via the process described in\\
\> \> \> Section \ref{sec:leaf}.\\
(5) \> \textbf{end parfor}\\
Box merging\\
(6)\> Set $N_\text{merged} = N_\text{boxes}^L$. \\
(7)\> Set $N_\text{unmerged} = N_\text{boxes} - N_\text{boxes}^L$. \\
(8) \>\textbf{for} $l = L - 1,\, L - 2,\, \dots,\, 1$  $\quad$ (loop over levels) \\
(9) \>\> Set outer threads to $\theta_o^l$. \\
(10) \> \>\textbf{parfor} $\tau = N_\text{unmerged} - N_\text{boxes}^l + 1, \, N_\text{unmerged} - N_\text{boxes}^l + 2, \, \dots, \, N_\text{unmerged}$ \\
(11) \> \> \> Set inner threads to $\theta_i^l$. \\
(12)\> \> \> Let $\alpha$ and $\beta$ be the children of $\tau$.\\
(13)\> \> \> Split $I_{\rm b}^{\alpha}$ and $I_{\rm b}^{\beta}$ into vectors $I_{1}$, $I_{2}$, and $I_{3}$ as shown in \\
\> \> \> Figure \ref{fig:siblings_notation}.\\
(14)\> \> \> $\mtx{W} = \mtx{I}-\mtx{R}^{\beta}_{33}\mtx{R}^{\alpha}_{33}$\\
(15)\> \> \> $\mtx{\Phi}^{\alpha} = \mtx{W}^{-1}\left[\mtx{R}_{33}^{\beta} \mtx{R}_{31}^{\alpha} | -\mtx{R}_{32}^{\beta}\right]$\\
(16) \> \> \> $\mtx{\Upsilon}^{\alpha} = \left[\mtx{W}^{-1}\mtx{R}_{33}^{\beta}\ | \ -\mtx{W}^{-1}\right] $\\
(17)\> \> \> $\mtx{\Phi}^{\beta} = \left[-\mtx{R}_{31}^{\alpha} -\mtx{R}^{\alpha}_{33}\mtx{W}^{-1}\mtx{R}_{33}^{\beta}\mtx{R}^{\alpha}_{31} | \mtx{R}_{33}^{\alpha} \mtx{W}^{-1}\mtx{R}_{32}^{\beta}\right] $\\ 
(18)\> \> \> $\mtx{\Upsilon}^{\beta} = \left[-\left(\mtx{I}+\mtx{R}^{\alpha}_{33}\mtx{W}^{-1}\mtx{R}_{33}^{\beta}\right)  \ | \ \mtx{R}_{33}^{\alpha} \mtx{W}^{-1}\right] $\\  
(19)\> \> \> $\mtx{R}^\tau = \mtwo{\mtx{R}_{11}^{\alpha}}{\mtx{0}}{\mtx{0}}{\mtx{R}_{22}^{\beta}} + \vtwo{\mtx{R}^{\alpha}_{13}}{\mtx{R}^{\beta}_{23}} \vtwo{\mtx{\Phi}^{\alpha}}{\mtx{\Phi}^{\beta}}$\\
(20)\> \> \> $\mtx{\Gamma}^\tau = \vtwo{\mtx{R}^{\alpha}_{13}}{\mtx{R}^{\beta}_{23}} \vtwo{\mtx{\Upsilon}^{\alpha}}{\mtx{\Upsilon}^{\beta}}$\\
(21)\> \> \> Delete $\mtx{R}^{\alpha}$ and $\mtx{R}^{\beta}$.\\
(22)\> \> \textbf{end  parfor}\\
(23) \> \> $N_\text{merged} = N_\text{merged} + N_\text{boxes}^l$. \\
(24) \> \> $N_\text{unmerged} = N_\text{unmerged} - N_\text{boxes}^l$. \\
(25)\> \textbf{end for}
\end{tabbing}
\end{minipage}}
\label{alg:parprecompute}
\end{figure}

\begin{figure}[hp]
\fbox{
\begin{minipage}{115mm}
\begin{center}
\textsc{Algorithm 4} (Solve stage, parallel)
\end{center}
This algorithm presents the parallel version 
of the solve stage.  The serial version is 
presented in Algorithm 2.
Let $N_\text{boxes}^l$ be 
the number of boxes on a level $l$, and $\theta_{ou}^l$, $\theta_{iu}^l$, $\theta_{od}^l$, and $\theta_{id}^l$ be the outer and inner upward pass stage threads and outer and inner downward pass stage threads 
for a level respectively.
All operators are assumed to be precomputed.

\rule{\textwidth}{0.5pt}
\begin{tabbing}
\mbox{}\hspace{7mm} \= \mbox{}\hspace{6mm} \= \mbox{}\hspace{6mm} \= \mbox{}\hspace{6mm} \= \mbox{}\hspace{6mm} \= \kill
Upward pass\\
(1) \> Set outer threads to  $\theta_{ou}^L$. \\
(2) \> \textbf{parfor} $\tau = N_\text{boxes} - N_\text{boxes}^L + 1,  \, N_\text{boxes} - N_\text{boxes}^L + 2, \dots, \, N_\text{boxes}$ \\
(3) \>\> Set inner threads to $\theta_{iu}^L$.\\
(4)\> \> Compute $\prtsolnvec^\tau = \PrtSolnOp^\tau\vct{s}$ and $\outimpprtvec^\tau = \OutImpPrtSolnOp^\tau\vct{s}$ for the leaf.\\
(5) \> \textbf{end parfor} \\
(6)\> Set $N_\text{processed} = N_\text{boxes}^L$. \\
(7)\> Set $N_\text{unprocessed} = N_\text{boxes} - N_\text{boxes}^L$. \\
(8) \>\textbf{for} $l = L - 1,\, L - 2,\, \dots,\, 1$  $\quad$ (loop upward through levels) \\
(9) \>\> Set outer threads to $\theta_{ou}^l$. \\
(10) \> \>\textbf{parfor} $\tau = N_\text{unprocessed} - N_\text{boxes}^l + 1, \, N_\text{unprocessed} - N_\text{boxes}^l + 2, \, \dots, \, N_\text{unprocessed}$ \\
(11) \> \> \> Set inner threads to $\theta_{iu}^l$. \\
(12)\> \> \> Let $\alpha$ and $\beta$ be the children of $\tau$.\\
(13)\> \> \> Compute $\tilde{\vct{t}}^{\alpha}_3 = \mtx{\Upsilon}^{\alpha}\vtwo{\vct{h}^{\alpha}_3}{\vct{h}^{\beta}_3}$ and $\tilde{\vct{t}}^{\beta}_3 = \mtx{\Upsilon}^{\beta}\vtwo{\vct{h}^{\alpha}_3}{\vct{h}^{\beta}_3}$.\\
(14)\> \> \> Compute $\outimpprtvec^\tau$ via (\ref{eq:flux_part}).\\
(15)\> \> \textbf{end parfor}\\
(16) \> \> $N_\text{processed} = N_\text{processed} + N_\text{boxes}^l$. \\
(17) \> \> $N_\text{unprocessed} = N_\text{unprocessed} - N_\text{boxes}^l$. \\
(18)\> \textbf{end for}\\
Downward pass\\
(19) \> Set $N_\text{processed} = 0$. \\
(20) \> Set $N_\text{unprocessed} = N_\text{boxes}$. \\
(21) \> \textbf{for} $l = 1, \, 2, \, \dots, \, L - 1$ $\quad$ (loop downward through levels) \\ 
(22) \> \> Set outer threads to $\theta_{od}^l$.\\
(23)\>\> \textbf{parfor} $\tau = N_\text{processed} + 1,\, N_\text{processed} + 2,\,\dots, \, N_\text{processed} + N_\text{boxes}^l$ \\
(24) \> \>\> Set inner threads to $\theta_{id}^l$. \\
(25)\> \> \> Let $\alpha$ and $\beta$ be the children of $\tau$.\\
(26)\> \> \>$\vct{t}^{\alpha}_3 = \IncImpHmgSolnOp^{\alpha}\vct{t}^\tau + \tilde{\vct{t}}^{\alpha}$, $\vct{t}^{\beta}_3 = \IncImpHmgSolnOp^{\beta}\vct{t}^\tau + \tilde{\vct{t}}^{\beta}$.\\
(27)\>\>\textbf{end parfor}  \\
(28) \>\> $N_\text{processed} = N_\text{processed} + N_\text{boxes}^l$.\\
(29) \> \> $N_\text{unprocessed} = N_\text{unprocessed} - N_\text{boxes}^l$. \\
(30) \>\textbf{end for} \\
(31) \> Set outer threads to $\theta_{od}^L$. \\
(32) \> \textbf{parfor} $\tau = N_\text{processed} + 1, \, N_\text{processed} + 2, \, \dots, \, N_\text{processed} + N_\text{boxes}^L$ \\
(33) \> \> Set inner threads to $\theta_{id}^L$. \\
(34)\>  \> Let $J^\tau$ denote the indices for the 
discretization points in $\tau$.  \\
(35)\> \>$\uu(J^\tau) =\HmgSolnOp^{\tau}\,\vct{t}^\tau + \tilde{\uu}^\tau$.\\
(36) \> \textbf{end parfor}
\end{tabbing}
\end{minipage}}
\label{alg:parsolve}
\end{figure}

We chose the number of inner and outer threads on a level  
based on the most expensive operation in processing
a box, called the \textit{representative action}.  In the build stage, the 
representative action is inverting a matrix.  For leaf boxes, this corresponds to the 
inverting the approximate differential operator.  For merging two boxes,
the inverse of the matrix defined on the interface is the representative
action.  In the solve stage, the representative action in both sweeps 
of the hierarchical tree is matrix-vector multiplication (\textit{matvec}).
Table \ref{tab:repactions} details the matrix size of the representative 
action for each stage based on level in the tree. 
We call the time for computing the representative action on level $l$
the \textit{representative time}.  Since this time depends on the 
number of threads given to parallel linear algebra, we 
denote the representative time for level $l$ with $j$ inner threads
by $r_l^j$.  Since the representative times are machine dependent
and the order of discretization order $n_c$ can be fixed for a 
variety of problems, the representative times may often only need to 
be found once for a machine.  The representative times
can then be used in the 
application of the HPS method to any boundary value problem of the form (\ref{eq:basic}).

\begin{center}
\begin{table}[h!]
\caption{Representative actions used to estimate computation time on a level.    }
\label{tab:repactions}
\begin{tabular}{c | c c c}

 Stage &   Level type           & matrix size & task \\ 
   \toprule
\multirow{2}{*}{Build} &  Leaf & interior $\times$ interior & \multirow{2}{*}{matrix inversion}  \\
 & Others &  interface $\times$ interface & \\
\midrule
\multirow{2}{*}{Upward pass} &  Leaf & [interior $+$ exterior] $\times$ interior  & \multirow{2}{*}{matvec}  \\
 &  Others &     exterior $\times$ interface &  \\
  \midrule
\multirow{2}{*}{Downward pass} &   Leaf  & [interior $+$ exterior] $\times$ exterior &   \multirow{2}{*}{matvec} \\
  &  Others &    interface $\times$ exterior  & \\
\end{tabular}
\end{table}
\end{center}

%For any level $l$, the set of feasible combinations of outer $\theta_o$ 
%and inner $\theta_i$ threads, $S^l$, satisfy 
% The possible combinations of $\theta_o$ and $\theta_i$ are chosen such that 
%\[ \theta_o\theta_i \leq \theta_t, \ \rm{and} \   \theta_o, \theta_i \in \mathbb{Z}. \]

For a uniform discretization the number of boxes on a given level is 
$N_\text{boxes}^l = 2^l$.  For a given
representative action, the choice
of inner and outer thread pair is selected so that 

% The estimated time $\epsilon_l$ on a level $l$ with $N_l$ boxes 
% for each possible thread combination is then
\begin{equation*}
 \epsilon^l=  \verb|ceiling|(N_\text{boxes}^l/\theta^l_o) \, r_l^{ \theta^l_i}
  \label{eq:reptime}
\end{equation*}
is minimized over all feasible pairs of $\theta^l_o$ and $\theta^l_i$.  Figures \ref{fig:opt-build-threads} and~\ref{fig:opt-solve-inner-threads} report the 
optimial choice of inner and outer threads for the 
build and solve stages of the algorithm with a total
of $\theta_t = 56$ total available threads on several
different uniform tree.  
 (The number of Chebyshev discretization points in each direction for  a leaf, $n_c$, is set at 16.  )

For the build stage, when $N_\text{boxes}^l > 10$, 
it is advantageous to let $\theta^l_i = 1$ and to distribute 
the boxes among the threads.  At the top several 
levels of the tree when the number of boxes on the
level is less than or equal to $10$, it is advantageous 
to assign each thread a box, i.e. $\theta^l_o = N_\text{boxes}^l$, and 
divide the remaining threads evenly for parallel linear algebra.  
Linear algebra is not parallelized in the lower levels of the tree
since the matrices are too small to benefit from it.  A similar behavior is observed for the solve stage, though
the parallel linear algebra is utilized earlier in the 
hierarchical tree, when there are approximately fewer than 15 
boxes on a level. 
\begin{figure}[h!]
     \begin{subfigure}{0.48\textwidth}
     \def\picsize{100}
     \def\xlabelloc{40}
     \def\ylabelloc{20}
     \begin{picture}(\picsize, \picsize)
       \put(10,10){\includegraphics[width=0.9\textwidth]{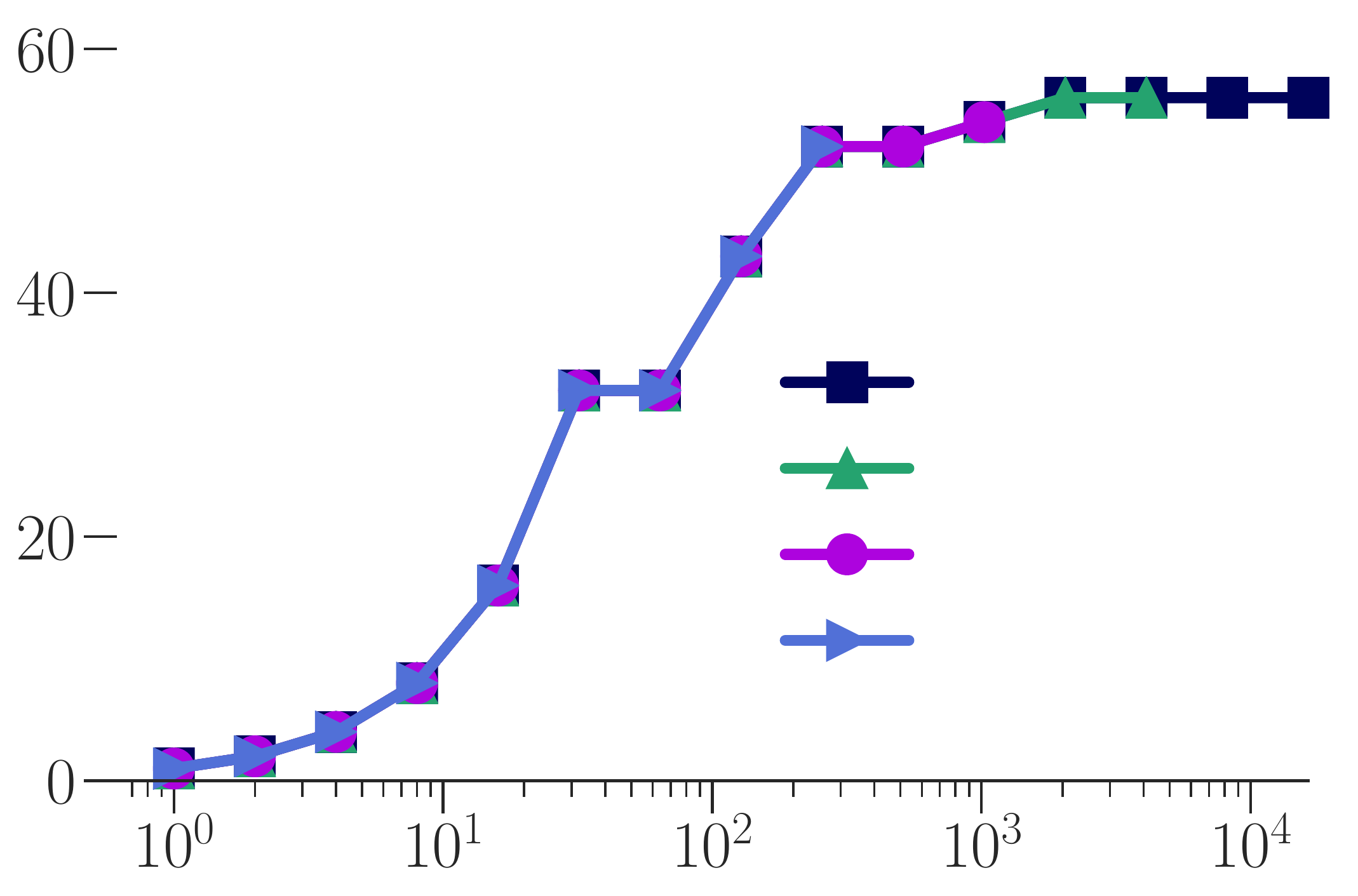}}
       \put(\xlabelloc, 0){number of boxes on level}
       \put(0, \ylabelloc){\rotatebox{90}{number of threads}}
       \put(112,64){\scriptsize{15 \textit{total levels}}}
        \put(112,55){\scriptsize{13}}
         \put(112,46){\scriptsize{11}}
         \put(112,37){\scriptsize{9}}
       \end{picture}
       \caption{\hspace{1cm}}
     \end{subfigure}\hfill
     \begin{subfigure}{0.48\textwidth}
      \def\picsize{100}
     \def\xlabelloc{40}
     \def\ylabelloc{20}
     \begin{picture}(\picsize, \picsize)
       \put(10,10){ \includegraphics[width=0.9\textwidth]{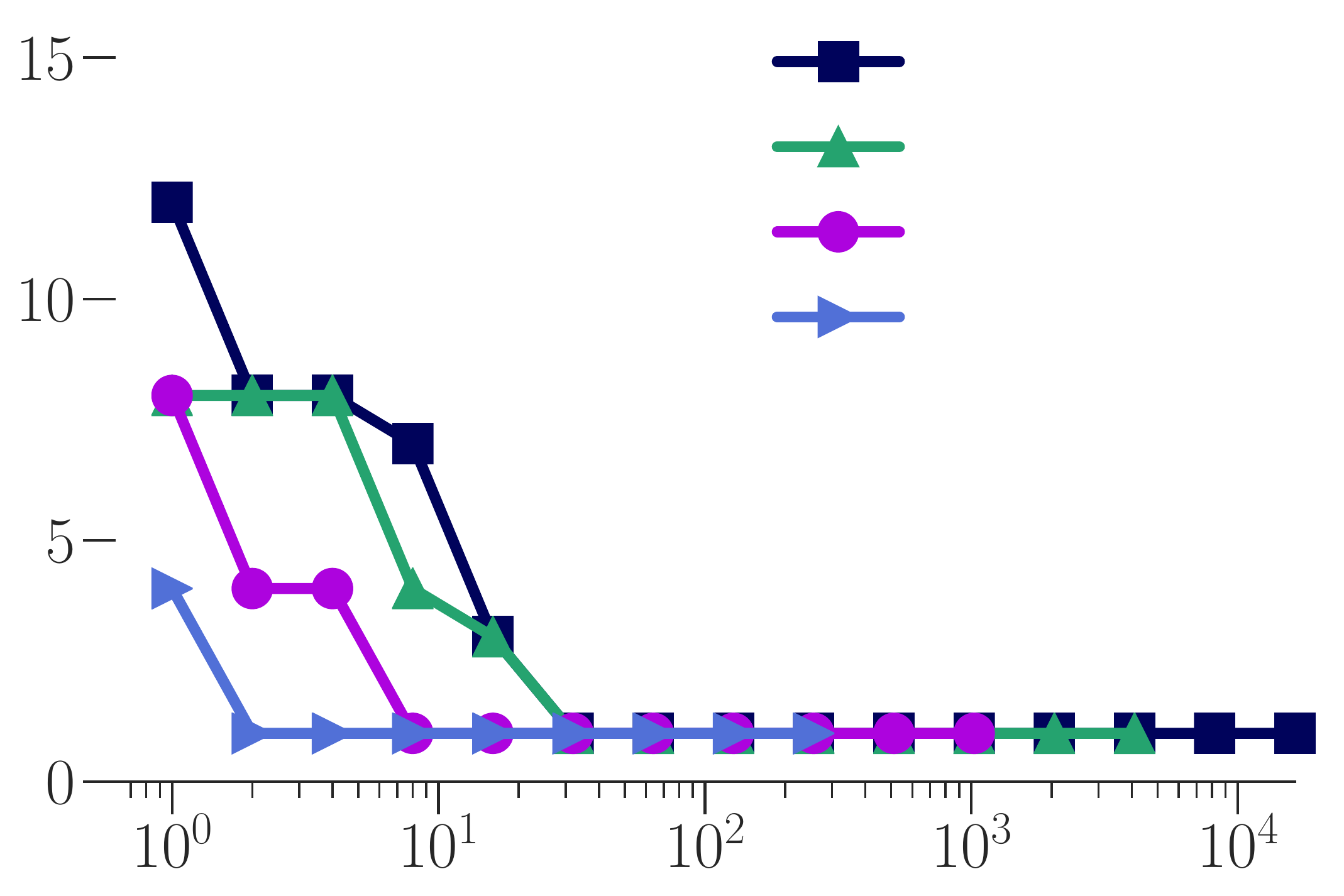}}
        \put(\xlabelloc, 0){number of boxes on level}
       \put(0, \ylabelloc){\rotatebox{90}{number of threads}}
         \put(115,100){\scriptsize{15 \textit{total levels}}}
        \put(115,91){\scriptsize{13}}
         \put(115,82){\scriptsize{11}}
         \put(115,73){\scriptsize{9}}
       \end{picture}
       \caption{\hspace{1cm}}
     \end{subfigure}
      \caption{The optimal combinations per level of (a) outer and (b) inner threads
      in the build stage of the HPS method with $9$ ({\color{lightblue}$\blacktriangleright$}), 
      $11$ ({\color{magenta} $\bullet$}), $13$ (
      {\color{green} $\blacktriangle$}), and $15$ ({\color{darkblue}$\blacksquare$})
      total levels in the tree, with a leaf discretization order of $n_c = 16$ and $56$ available threads.}
  \label{fig:opt-build-threads}
   \end{figure}
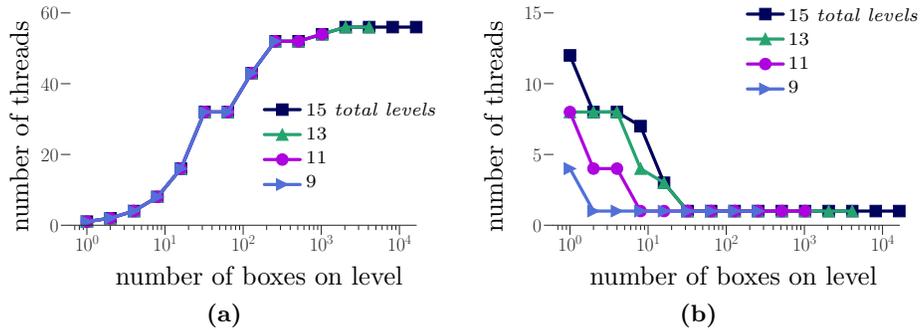
   \begin{figure}[h!]
     \begin{subfigure}{0.48\textwidth}
     \def\picsize{100}
     \def\xlabelloc{40}
     \def\ylabelloc{20}
     \begin{picture}(\picsize, \picsize)
       \put(10,10){ \includegraphics[width=0.9\textwidth]{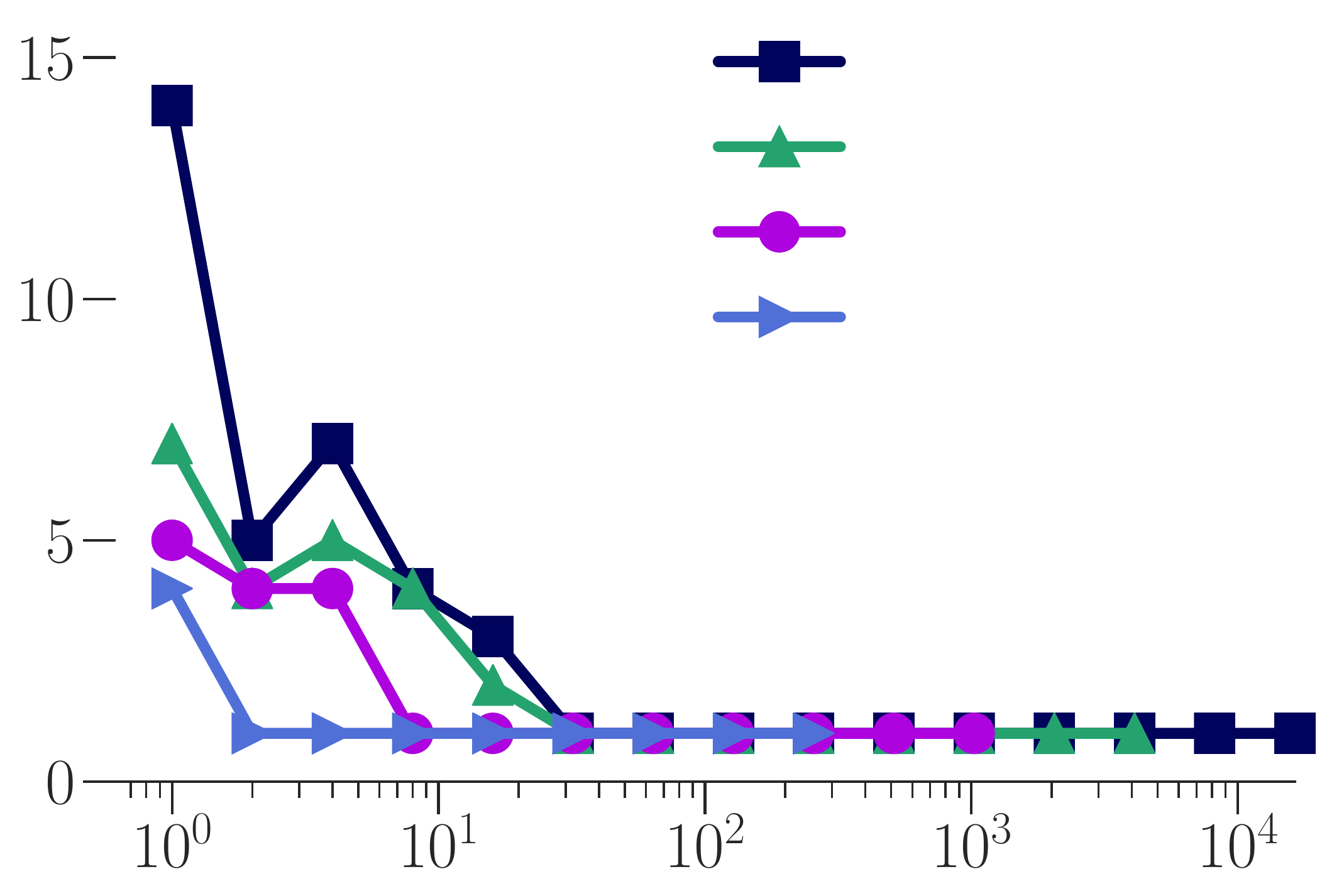}}
        \put(\xlabelloc, 0){number of boxes on level}
       \put(0, \ylabelloc){\rotatebox{90}{number of threads}}
          \put(112,100){\scriptsize{15 \textit{total levels}}}
        \put(112,91){\scriptsize{13}}
         \put(112,82){\scriptsize{11}}
         \put(112,73){\scriptsize{9}}
       \end{picture}
       \caption{}
        \label{fig:solve-uw-opt-inner-threads}
     \end{subfigure}\hfill
     \begin{subfigure}{0.48\textwidth}
      \def\picsize{100}
     \def\xlabelloc{40}
     \def\ylabelloc{20}
     \begin{picture}(\picsize, \picsize)
       \put(10,10){ \includegraphics[width=0.9\textwidth]{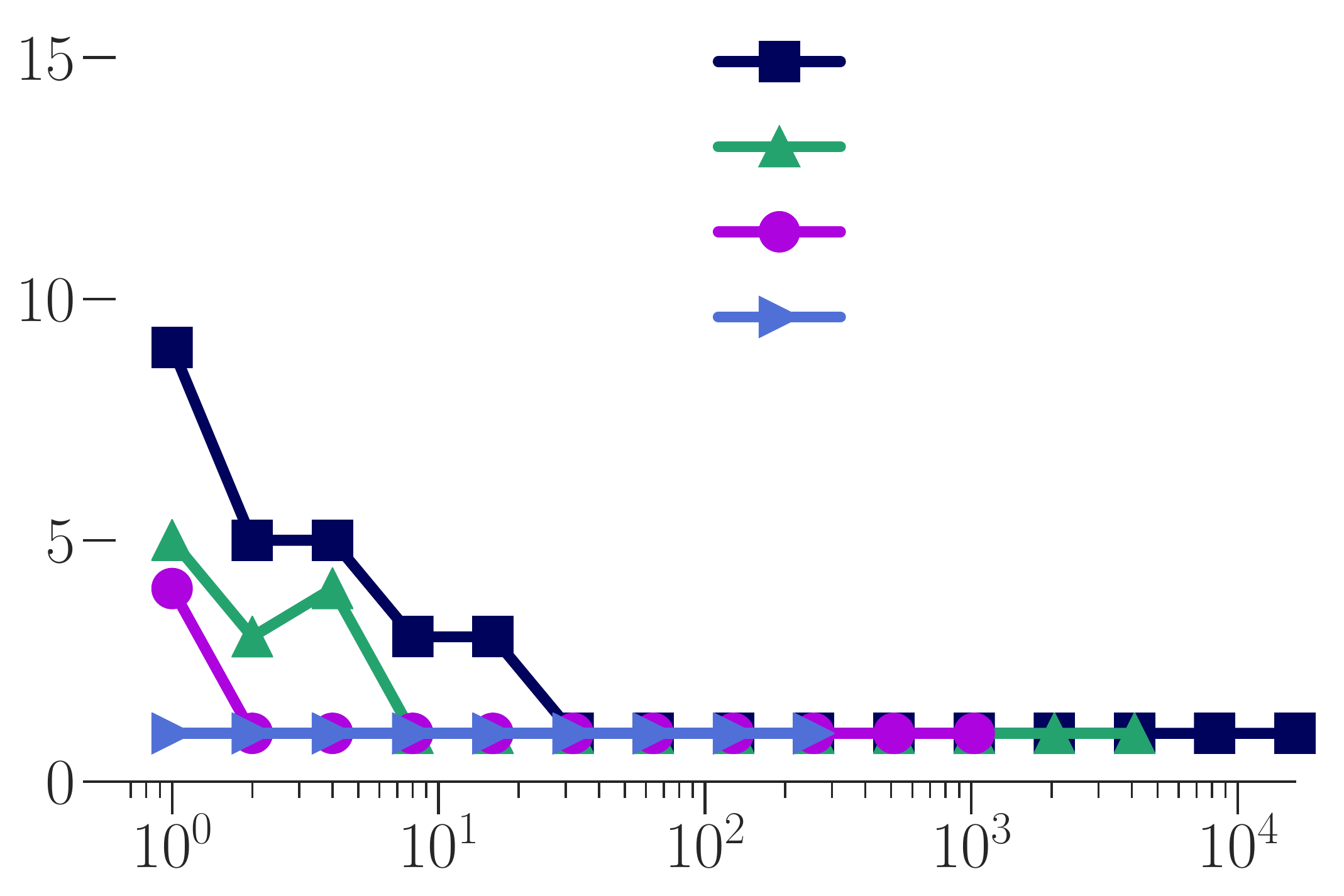}}
     \put(\xlabelloc, 0){number of boxes on level}
       \put(0, \ylabelloc){\rotatebox{90}{number of threads}}
        \put(112,100){\scriptsize{15 \textit{total levels}}}
        \put(112,91){\scriptsize{13}}
         \put(112,82){\scriptsize{11}}
         \put(112,73){\scriptsize{9}}
       \end{picture}
       \caption{}
        \label{fig:solve-dw-opt-inner-threads}
     \end{subfigure}
      \caption{Optimal inner threads on each level of the tree for upward (a) and downward (b) passes of the solve stage with 56 available threads.  Total tree depths are  $9$ ({\color{lightblue}$\blacktriangleright$}), 
      $11$ ({\color{magenta} $\bullet$}), $13$ (
      {\color{green} $\blacktriangle$}), and $15$ ({\color{darkblue}$\blacksquare$}), with a leaf discretization order of $n_c = 16$.  The optimal
      number of outer threads on each level is identical to the build stage, shown in Figure \ref{fig:opt-build-threads}(a).}
           \label{fig:opt-solve-inner-threads}
   \end{figure}
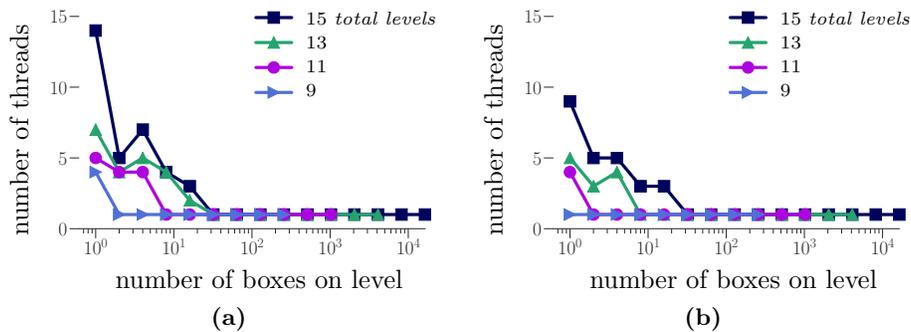
\begin{remark}
When the hierarchical tree in the HPS method is non-uniform, 
the size of the largest matrix on a level $l$ can be used 
to create the calibration data.  If the tree is highly
non-uniform and is going to remain fixed for a large number 
of simulations, it may be advantageous to build calibration data
for each possible matrix size and adjust the definition of 
$\epsilon^l$ accordingly.

\end{remark}

 \section{Results}
  \label{sec:results}

This section illustrates the performance of the 
parallelization technique when implemented on 
a dual 2.3 GHz Intel Xeon Processor E5-2695 v3 desktop workstation 
with 256 GB of RAM, 28 physical cores, and $56$ possible threads.  The 
algorithm was implemented in Fortran 95, with the gfortran 5.4.0 compiler
using the \verb|-O3| optimization flag.  MKL vectorization was not used.

The HPS method was applied to (\ref{eq:basic}) where $\Omega$ is the 
unit square, $\kappa = 16$, the coefficent function $c(\vct{x}) = \exp\{-8 (x-0.5)^2 +(y-0.5)^2\}$ is 
a Gaussian centered in $\Omega$ and the exact solution is    $u(\vct{x})=u(x,y) = \exp{(i2\pi \kappa x)}\exp{(i2\pi \kappa y)}$.  For all experiments, $N$ denotes the number of discretization 
   points where the solution is unknown. 

Sections \ref{sec:results-nc} and \ref{sec:results-nt} report the performance of the
parallel implementation with varying the discretization order (via $n_c$) and the number of available threads $\theta_t$, respectively.  Section \ref{sec:results-roofline} reports the hardware efficiency 
of the parallel implementation.

\subsection{Tests with varying $n_c$}
\label{sec:results-nc}
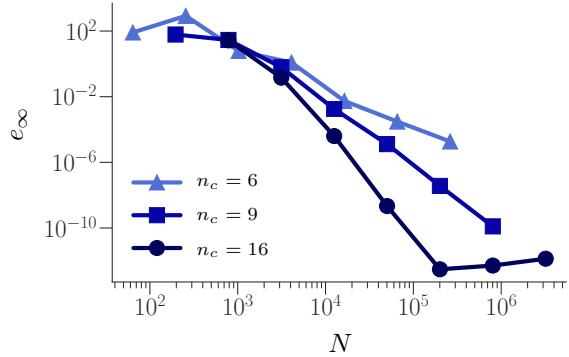
\begin{figure}[h!]
      \vspace*{0.2\textwidth}
        \begin{center} 
      \def\picsize{100}
     \def\xlabelloc{40}
     \def\ylabelloc{80}
     \begin{picture}(\picsize, \picsize)
       \put(-70,10){\includegraphics[width=0.6\textwidth]{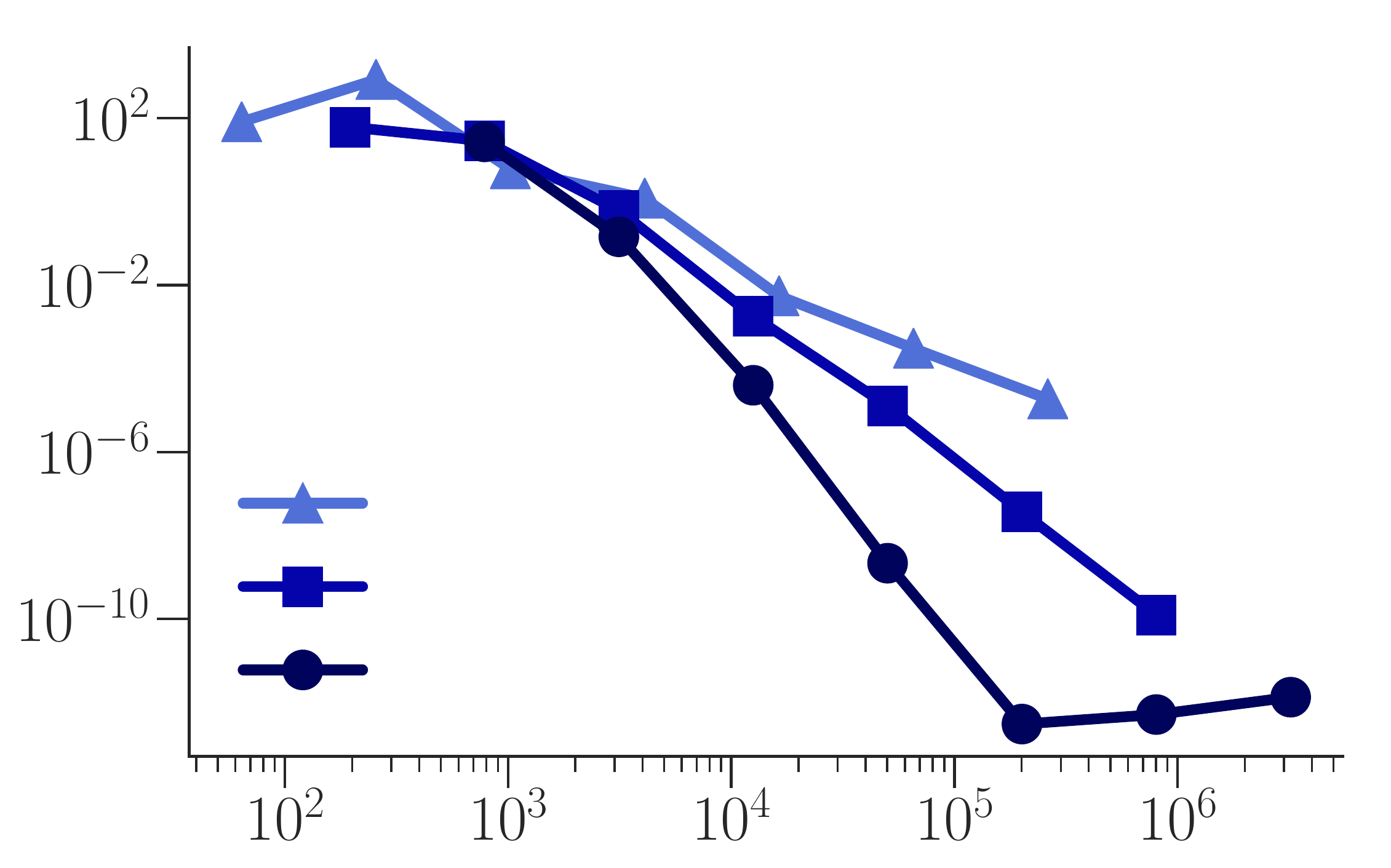}}
      \put(\xlabelloc, 0){$N$}
       \put(-80, \ylabelloc){\rotatebox{90}{$e_\infty$}}
            \put(-10,63){\scriptsize{$n_c = 6$}}
            \put(-10,50){\scriptsize{$n_c = 9$}}
            \put(-10,37){\scriptsize{$n_c = 16$}}
       \end{picture}
       \end{center}
      \caption{Maximum error $e_\infty$ in the approximate solution versus total number of interior unknowns $N$ for the same set of trees with varying order of Chebyshev discretization $n_c$ plotted on a log-log scale. }
           \label{fig:error}
      \end{figure}

In this section, three choices of order of discretization are considered; $n_c = 6,\ 9,$ and $16$.  Let $e_\infty = max_{j=1,\ldots,N}|u(\vct{x}_j)- \vct{u}(j)|$ where $\{\vct{x}_j\}_{j =1}^N$ are discretization points where the solution 
is not known and $\vct{u}$ is the $N\times 1$ vector whose $j^{\rm th}$
entry is the approximate solution at the discretization point $\vct{x}_j$.  
Figure \ref{fig:error} reports the error $e_\infty$ versus $N$
for each choice of $n_c$.  For similar number of unknowns $N$, $n_c = 16$ 
   achieves the best accuracy as $N$ is increased beyond the point where
   the discretization can begin to fully resolve the solution. 

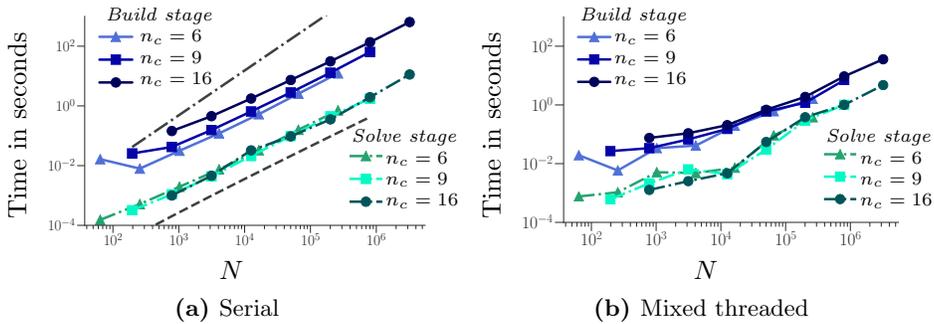
\begin{figure}[h!]
     \begin{subfigure}{0.48\textwidth}
       \def\picsize{100}
     \def\xlabelloc{80}
     \def\ylabelloc{23}
     \begin{picture}(\picsize, \picsize)
       \put(10,10){ \includegraphics[width=0.9\textwidth]{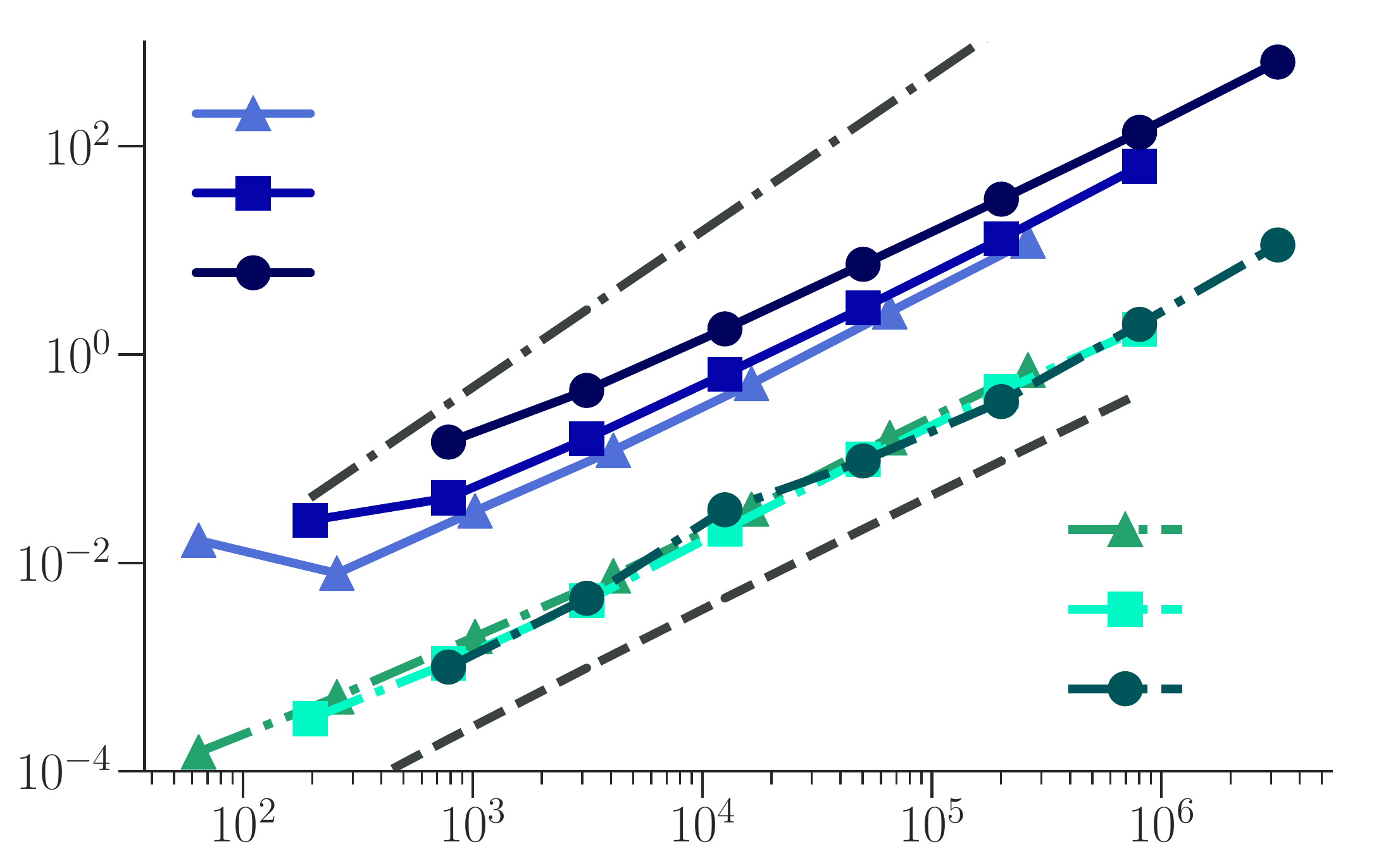}}
     \put(\xlabelloc, 0){$N$}
       \put(0, \ylabelloc){\rotatebox{90}{Time in seconds}}
       \put(37, 98){\scriptsize{\textit{Build stage}}}
       \put(50,90){\scriptsize{$n_c = 6$}}
        \put(50,82){\scriptsize{$n_c = 9$}}
        \put(50,74){\scriptsize{$n_c = 16$}}
        \put(130, 52){\scriptsize{\textit{Solve stage}}}
         \put(143,44){\scriptsize{$n_c = 6$}}
        \put(143,36){\scriptsize{$n_c = 9$}}
        \put(143,28){\scriptsize{$n_c = 16$}}
       \end{picture}
        \label{fig:serial-times}
       \caption{Serial}
     \end{subfigure}\hfill
     \begin{subfigure}{0.48\textwidth}
        \def\picsize{100}
     \def\xlabelloc{80}
     \def\ylabelloc{23}
     \begin{picture}(\picsize, \picsize)
       \put(10,10){ \includegraphics[width=0.9\textwidth]{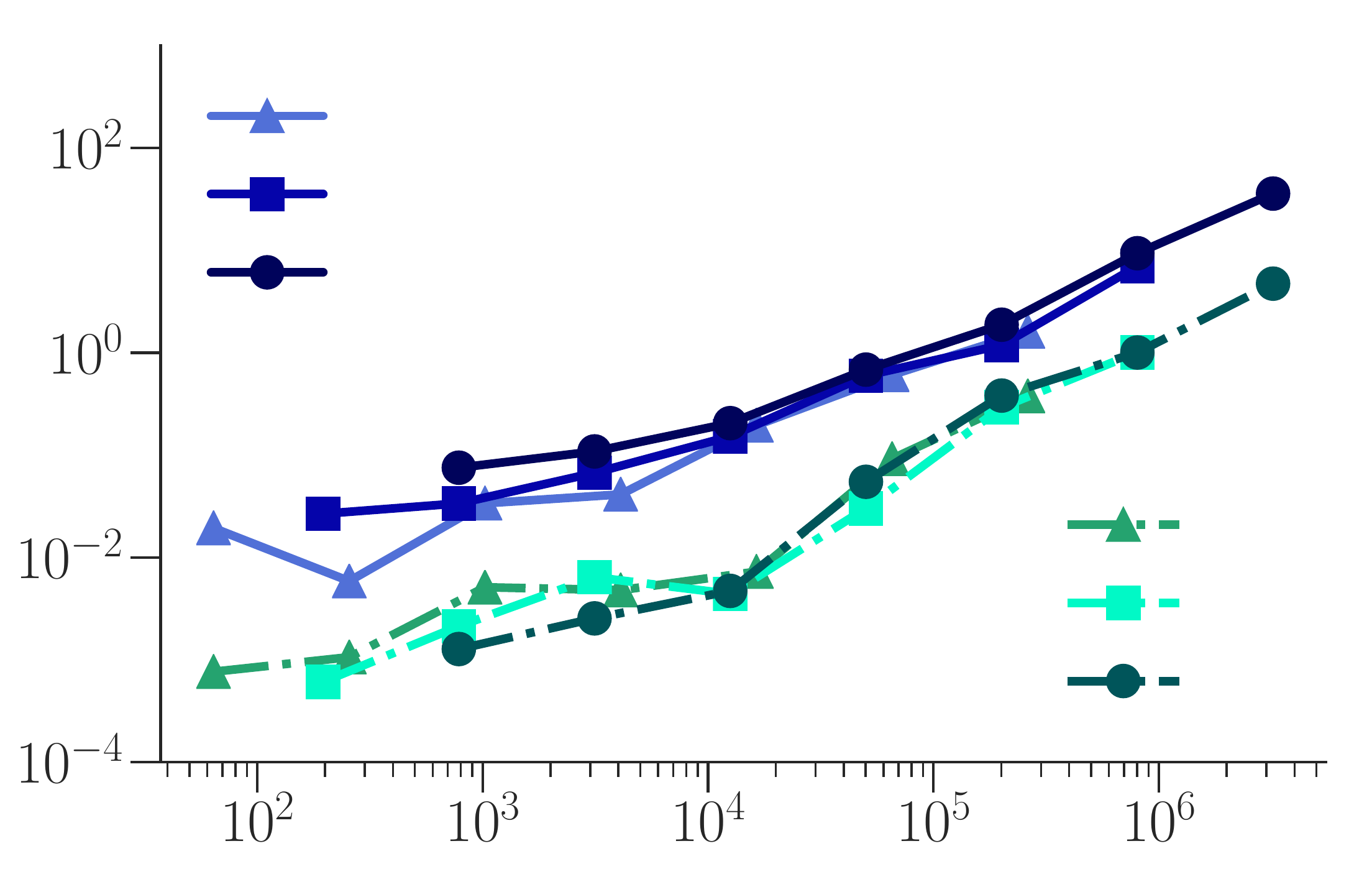}}
            \put(\xlabelloc, 0){$N$}
       \put(0, \ylabelloc){\rotatebox{90}{Time in seconds}}
              \put(37, 98){\scriptsize{\textit{Build stage}}}
         \put(50,90){\scriptsize{$n_c = 6$}}
        \put(50,82){\scriptsize{$n_c = 9$}}
        \put(50,74){\scriptsize{$n_c = 16$}}
               \put(130, 52){\scriptsize{\textit{Solve stage}}}
         \put(143,44){\scriptsize{$n_c = 6$}}
        \put(143,36){\scriptsize{$n_c = 9$}}
        \put(143,28){\scriptsize{$n_c = 16$}}
       \end{picture}
        \label{fig:threaded-times}
       \caption{Mixed threaded}
     \end{subfigure}
      \caption{Time in seconds versus the number of 
      discretization points where the solution is unknown $N$ for the (a) serial  and (b) mixed threaded ($\theta_t = 56$) 
      implementation for the precomputation and solve stages for varying values of $n_c$. Plotted on log-log scales.  In (a), the dashed-dot and dashed lines represent the $O(N^{1.5})$ and $O(N\log N)$ asymptotic complexities of the build and solve stages of the algorithm, respectively. }
      \label{fig:times}
   \end{figure}

 Figure  \ref{fig:times} reports 
   the execution time, in seconds, for the HPS method
   with the serial and 56-threaded parallel implementations. For a given problem size $N$, the parallel implementation nearly eliminates
   the additional cost in the build stage associated with higher $n_c$.
   Thus, the parallel implementation allows the HPS method to achieve higher accuracy at a very small cost in terms of additional computation time, as seen in Figure \ref{fig:error}.  This is because, in terms of $N$ and the ability to distribute work, order only impacts the size of the matrices at the leaf level.

\begin{figure}[h!]
     \begin{subfigure}{0.48\textwidth}
     \def\picsize{100}
     \def\xlabelloc{80}
     \def\ylabelloc{50}
     \begin{picture}(\picsize, \picsize)
       \put(10,10){\includegraphics[width=0.9\textwidth]{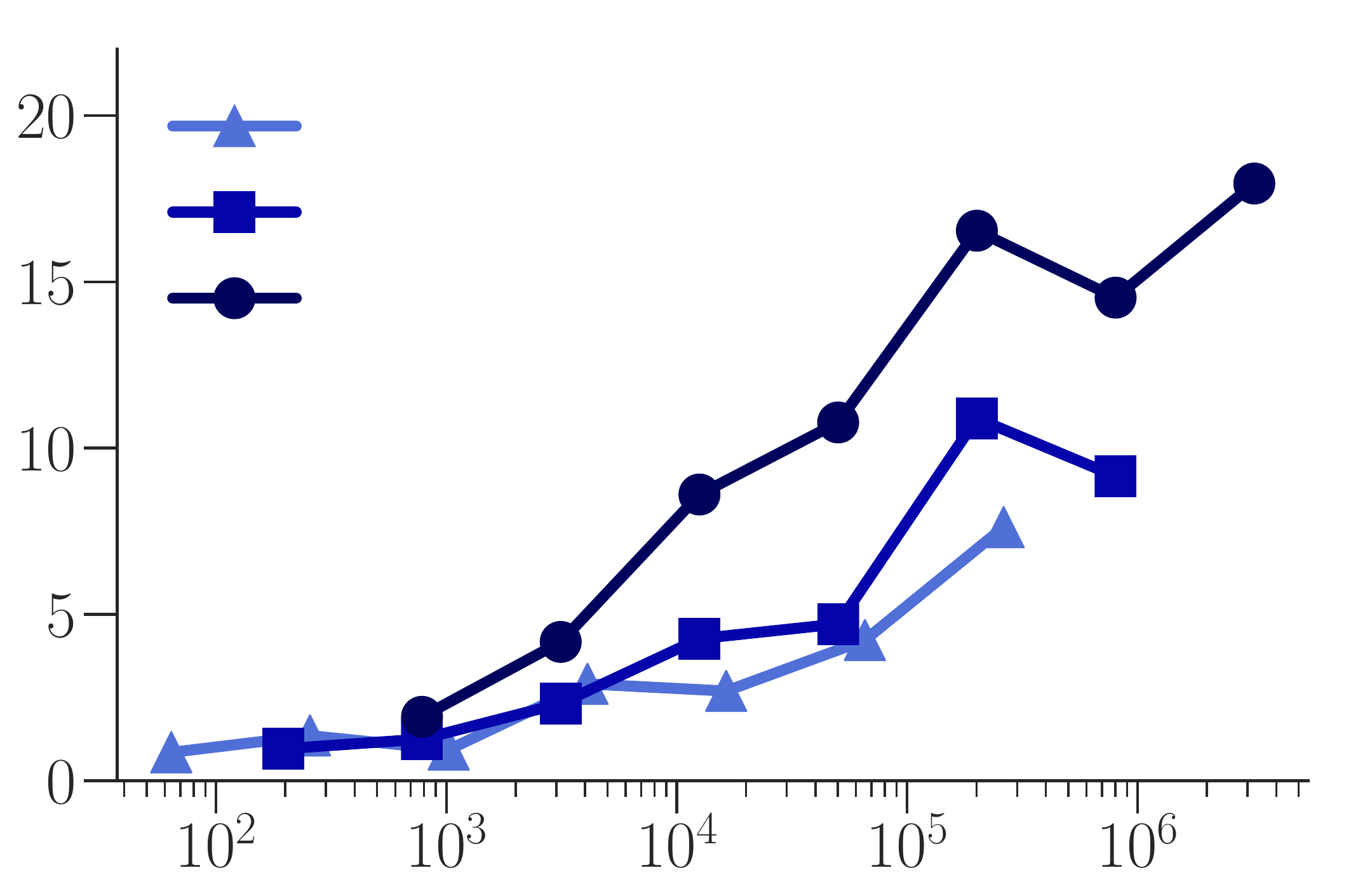}}
        \put(\xlabelloc, 0){$N$}
       \put(0, \ylabelloc){\rotatebox{90}{speedup}}
                \put(48,92){\scriptsize{$n_c = 6$}}
        \put(48,84){\scriptsize{$n_c = 9$}}
        \put(48,76){\scriptsize{$n_c = 16$}}
       \end{picture}
       \caption{Build stage}
     \end{subfigure}\hfill
     \begin{subfigure}{0.48\textwidth}
      \def\picsize{100}
     \def\xlabelloc{80}
     \def\ylabelloc{50}
     \begin{picture}(\picsize, \picsize)
       \put(10,10){\includegraphics[width=0.9\textwidth]{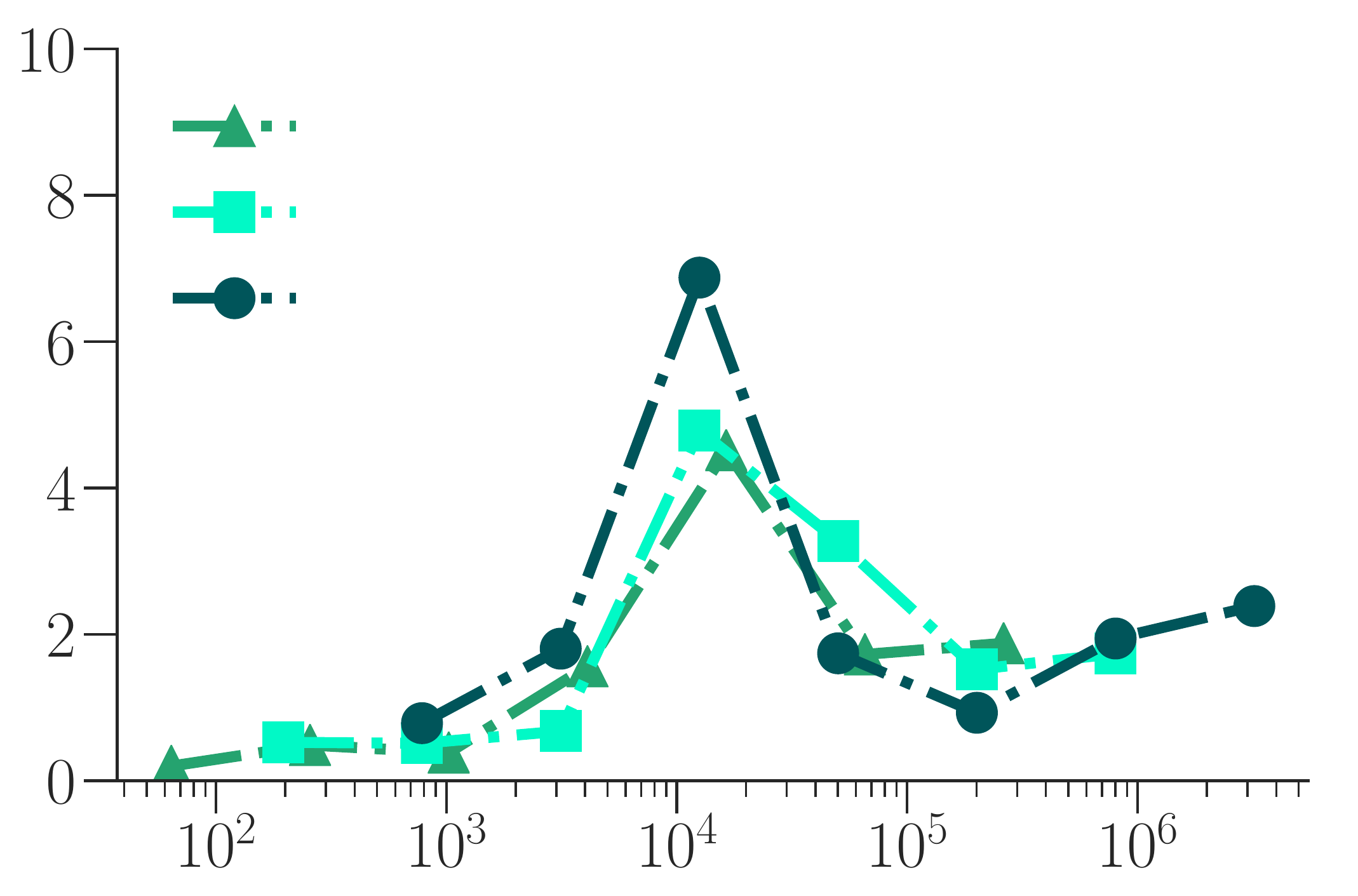}}
      \put(\xlabelloc, 0){$N$}
       \put(0, \ylabelloc){\rotatebox{90}{speedup}}
                       \put(48,92){\scriptsize{$n_c = 6$}}
        \put(48,84){\scriptsize{$n_c = 9$}}
        \put(48,76){\scriptsize{$n_c = 16$}}
       \end{picture}
       \caption{Solve stage}
     \end{subfigure}
      \caption{ Speedup in (a) build and (b) solve stage portions of the algorithm for 56 total threads with varying $n_c$ plotted on a semi-log scale.}
           \label{fig:speedup}
   \end{figure}
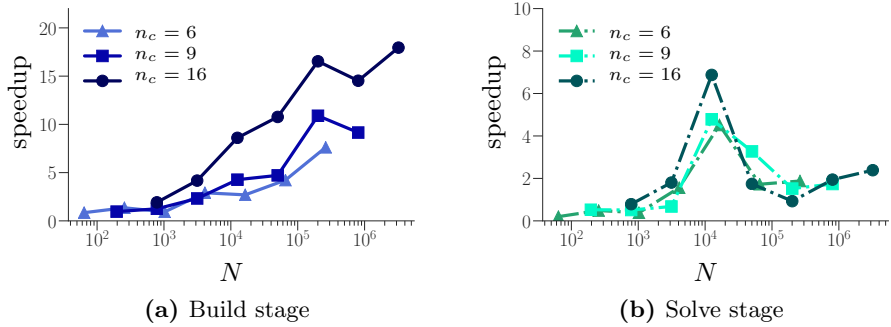

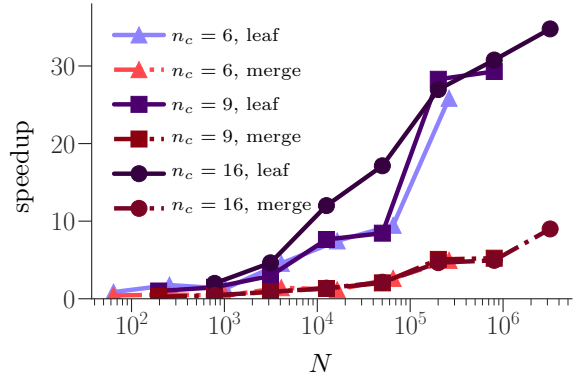
\begin{figure}[h!]
            \vspace*{0.05\textwidth}
            \hspace*{0.15\textwidth}
      \def\picsize{130}
     \def\xlabelloc{110}
     \def\ylabelloc{60}
     \begin{picture}(\picsize, \picsize)
       \put(10,10){\includegraphics[width=0.6\textwidth]{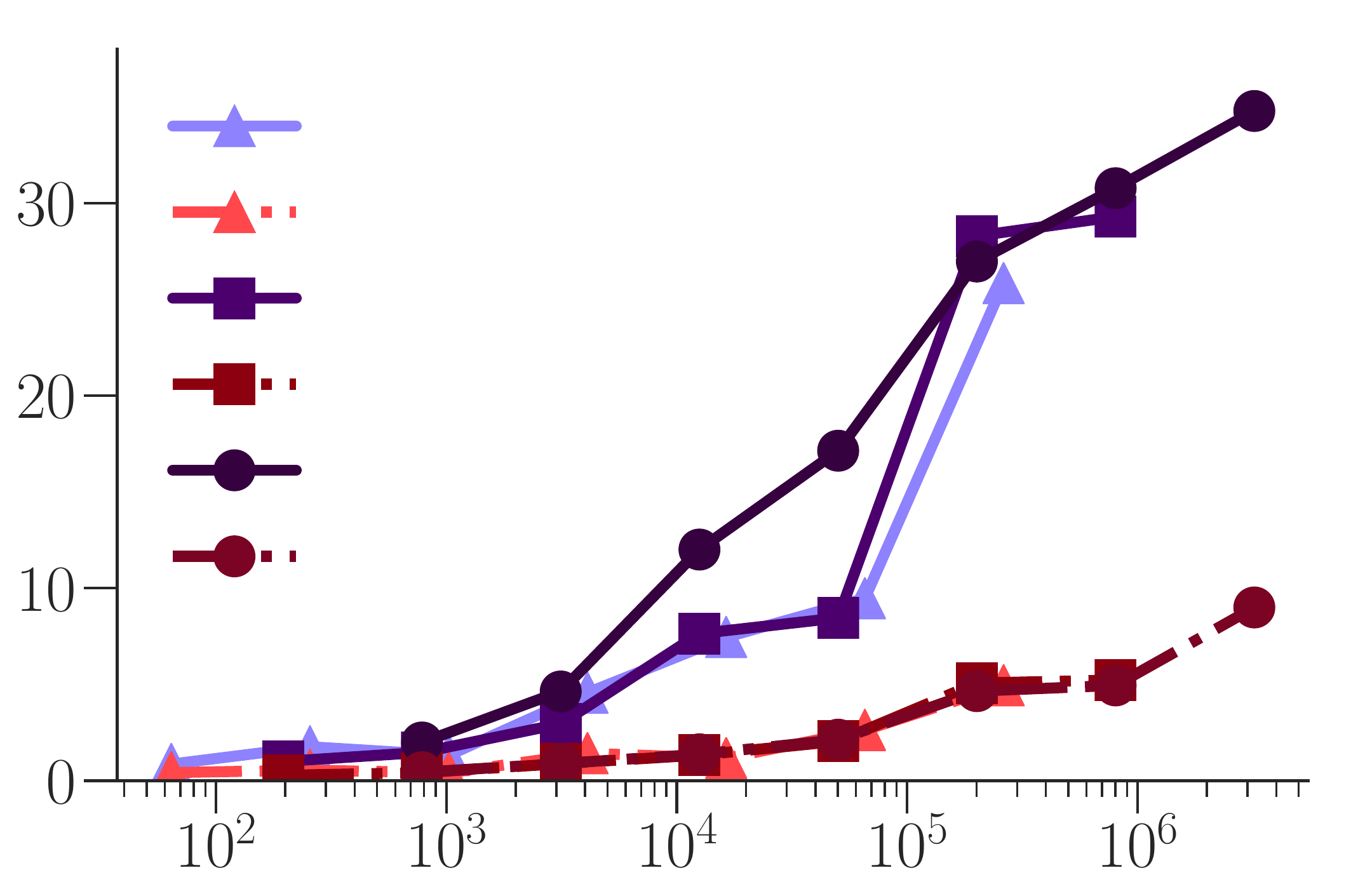}}
      \put(\xlabelloc, 0){$N$}
       \put(-3, \ylabelloc){\rotatebox{90}{speedup}}
            \put(58,125){\scriptsize{$n_c = 6$, leaf}}
            \put(58,112){\scriptsize{$n_c = 6$, merge}}    
        \put(58,99){\scriptsize{$n_c = 9$, leaf}}
        \put(58,87){\scriptsize{$n_c = 9$, merge}}
        \put(58,74){\scriptsize{$n_c = 16$, leaf}}
       \put(58,61){\scriptsize{$n_c = 16$,  merge}}
       \end{picture}
      \caption{Speedup in build stage of the algorithm, broken down by leaf and merge portions, for 56 total threads.}
           \label{fig:speedup-build-only}
      \end{figure}

  Figure \ref{fig:speedup}
   reports the corresponding speedup gained by moving from the serial
   to parallel implementation of the algorithm.  For the largest 
   problem considered with over two million unknowns, the build
   stage of the algorithm takes roughly 10 minutes via the 
   serial implementation while the parallel implementation takes 
   36 seconds, a parallel speedup of 17.5.    The constant prefactor for the solve stage is small since it is 
   a collection of BLAS3 matvec operations involving small matrices.  
   These efficient operations are precisely why the expected speedup is 
   small.  In fact, only an approximate factor of two speedup is obtained for the largest tree.
   Figure \ref{fig:speedup-build-only} reports the speedup, split between
   the leaf and merge computations in the build stage of the algorithm.
   As expected, the bulk of the speed up is gained on the leaf level
   where the algorithm is perfectly parallelizable.  The speedup
   of the merge computations is limited by the speedup of the
   parallel linear algebra inversion provided by MKL.     
  
   \subsection{Tests with varying $\theta_t$}
   \label{sec:results-nt}
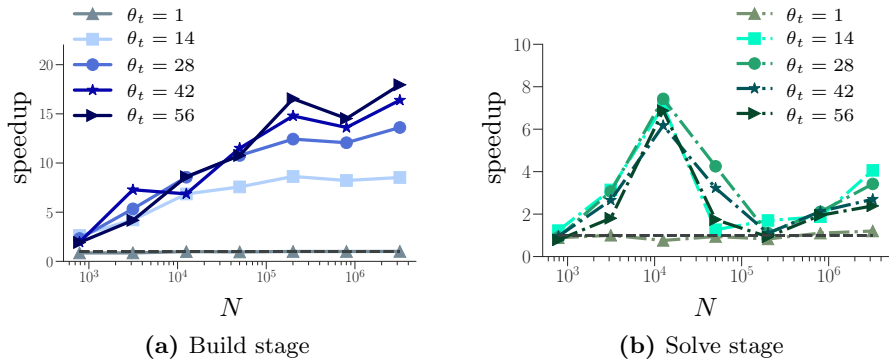
\begin{figure}[h!]
\vspace*{0.05\textwidth}
     \begin{subfigure}{0.48\textwidth}
     \def\picsize{100}
     \def\xlabelloc{80}
     \def\ylabelloc{50}
     \begin{picture}(\picsize, \picsize)
       \put(10,10){\includegraphics[width=0.9\textwidth]{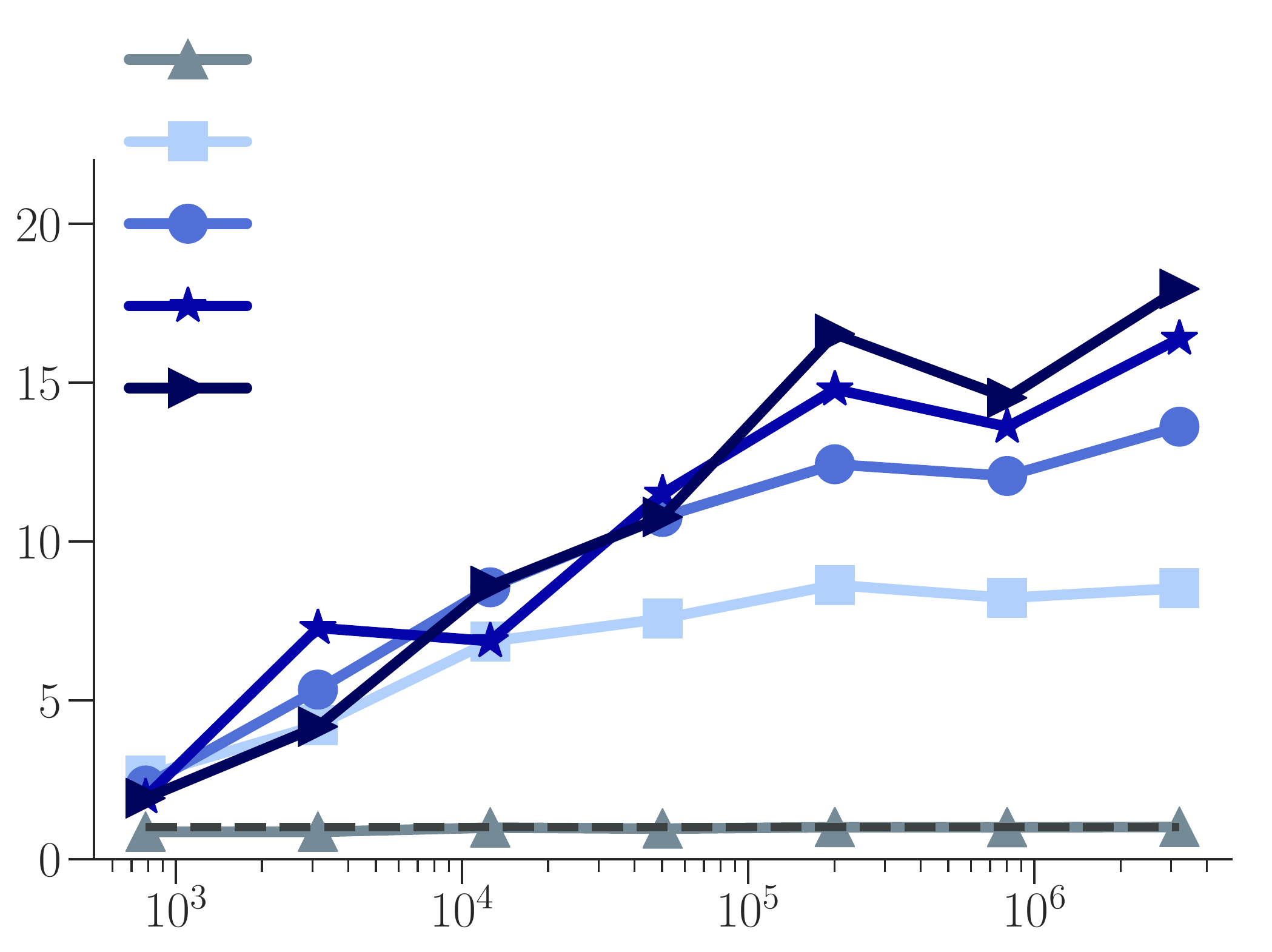}}
        \put(\xlabelloc, 0){$N$}
       \put(0, \ylabelloc){\rotatebox{90}{speedup}}
                \put(45,112){\scriptsize{$\theta_t = 1$}}
        \put(45,103){\scriptsize{$\theta_t = 14$}}
        \put(45,93){\scriptsize{$\theta_t = 28$}}
       \put(45,83){\scriptsize{$\theta_t = 42$}}
       \put(45,74){\scriptsize{$\theta_t = 56$}}
       \end{picture}
       \caption{Build stage}
     \end{subfigure}\hfill
     \begin{subfigure}{0.48\textwidth}
      \def\picsize{100}
     \def\xlabelloc{80}
     \def\ylabelloc{50}
     \begin{picture}(\picsize, \picsize)
       \put(10,10){\includegraphics[width=0.9\textwidth]{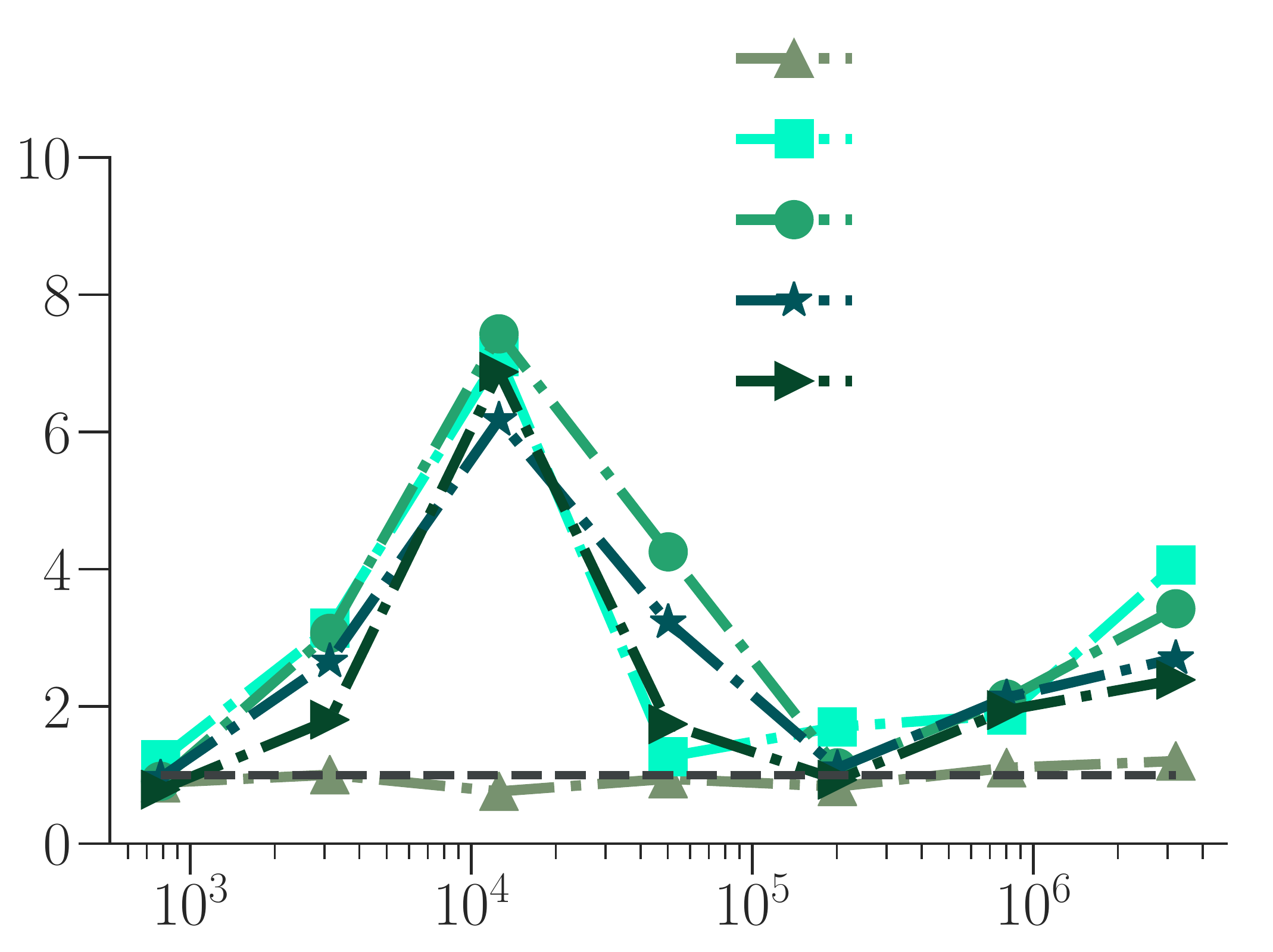}}
      \put(\xlabelloc, 0){$N$}
       \put(0, \ylabelloc){\rotatebox{90}{speedup}}
        \put(115,112){\scriptsize{$\theta_t = 1$}}
        \put(115,103){\scriptsize{$\theta_t = 14$}}
        \put(115,93){\scriptsize{$\theta_t = 28$}}
       \put(115,83){\scriptsize{$\theta_t = 42$}}
       \put(115,74){\scriptsize{$\theta_t = 56$}}
       \end{picture}
       \caption{Solve stage}
     \end{subfigure}
      \caption{Semi-log plots illustrating the parallel speedup in build (a) and solve (b) stage portions of the algorithm for varying number of total threads $\theta_t$ for
      $n_c = 16$. In both plots, the dashed lines denote speedup = 1, or equal time with the serial implementation.}
           \label{fig:speedup-var-threads}
   \end{figure}

   This section investigates the performance of the 
   parallel implementation for a varying number of 
   total threads $\theta_t$.  In these tests, 
   $n_c$ is fixed at $16$.
  Figure \ref{fig:speedup-var-threads} reports the speedup in the build 
  and solve stages for $\theta_t = \{1, 14, 28, 42, 56\}$. As expected,  increasing the number of threads
   improves build stage speedup for large problems; i.e. $N$ big.  
  For small problems, there are limited
  speed up gains as there is not much work to distribute amongst the threads even from the 
  divide-and-conquer strategy at the lower levels.
  At the upper levels in the tree, the matrices are modest sized and thus only experience modest speed up from threaded MKL linear algebra.  For larger problems, there is plenty of work at all levels in the tree.  There are performance gains from 
  the divide-and-conquer parallelism strategy on the lower levels of the tree and the matrices at the upper levels are large enough to benefit from threaded MKL linear algebra.    
  In the solve stage, a clear trend across the range of threads is less obvious, though the best performance
  for the largest problem tested is with just 14 threads.  This is due to the diminishing returns from increasing the number of threads for matvecs of the sizes required by the algorithm.  It is important to keep in mind that even the serial implementation of the solve stage takes less than a second of wall time for nearly all problem sizes tested.
      
   \subsection{Hardware efficiency}
   \label{sec:results-roofline}
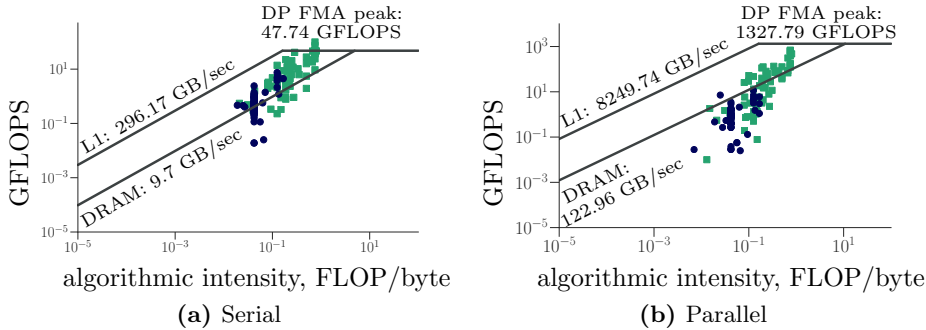
\begin{figure}[h!]
     \begin{subfigure}{0.48\textwidth}
     \def\picsize{100}
     \def\xlabelloc{20}
     \def\ylabelloc{30}
     \begin{picture}(\picsize, \picsize)
       \put(10,10){\includegraphics[width=0.9\textwidth]{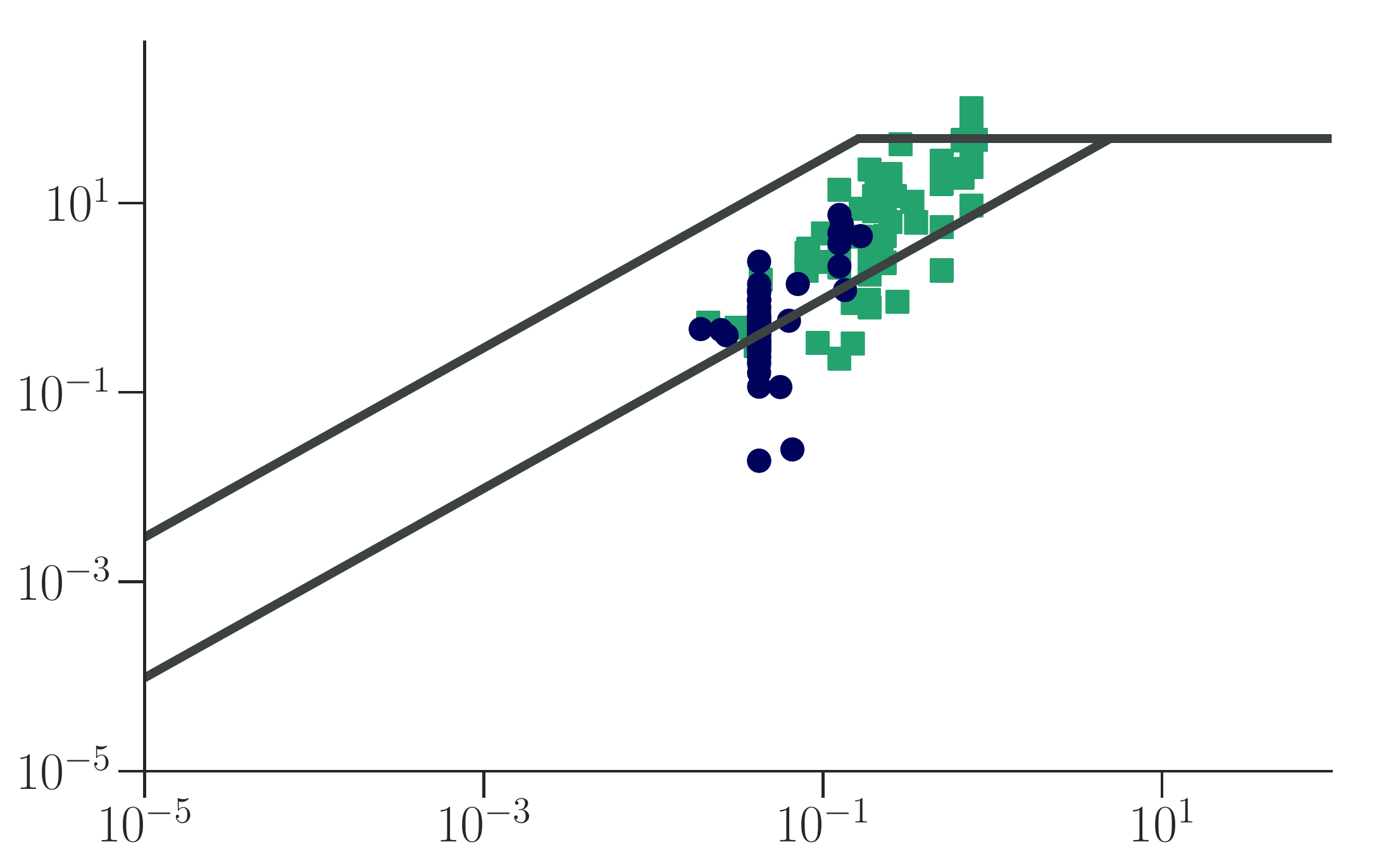}}
        \put(\xlabelloc, 0){algorithmic intensity, FLOP/byte}
       \put(0, \ylabelloc){\rotatebox{90}{GFLOPS}}
       \put(26, 50){\rotatebox{30}{\scriptsize{L1: 296.17 GB/sec}}}
        \put(26, 22){\rotatebox{30}{\scriptsize{DRAM: 9.7 GB/sec}}}
        \put(95,101){\scriptsize{DP FMA peak:}}
        \put(95,94){\scriptsize{47.74 GFLOPS}}
       \end{picture}
       \caption{Serial}
     \end{subfigure}\hfill
     \begin{subfigure}{0.48\textwidth}
      \def\picsize{100}
     \def\xlabelloc{20}
     \def\ylabelloc{30}
     \begin{picture}(\picsize, \picsize)
       \put(10,10){\includegraphics[width=0.9\textwidth]{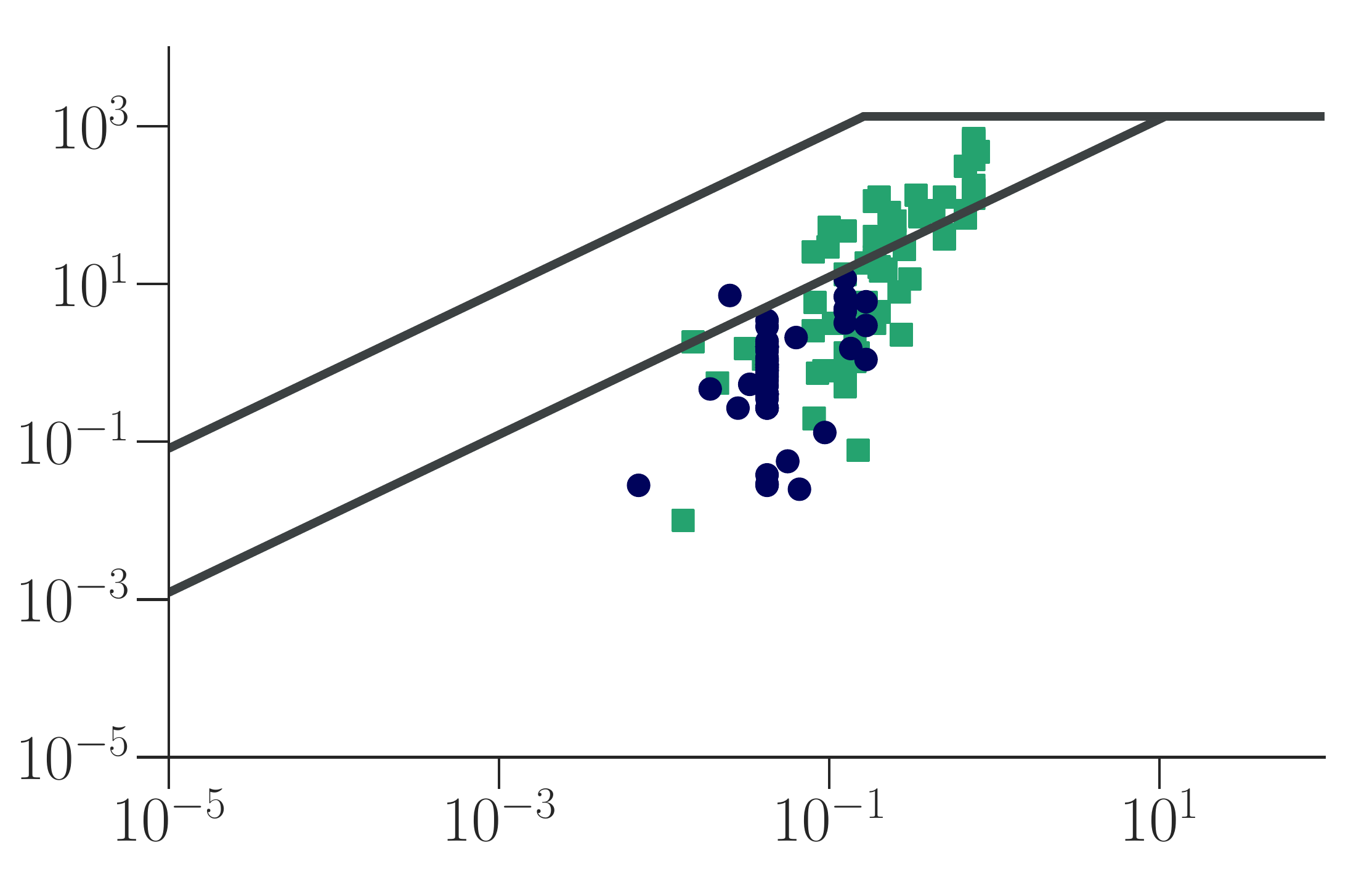}}
          \put(\xlabelloc, 0){algorithmic intensity, FLOP/byte}
       \put(0, \ylabelloc){\rotatebox{90}{GFLOPS}}
        \put(28, 60){\rotatebox{26}{\scriptsize{L1: 8249.74 GB/sec}}}
        \put(29, 31){\rotatebox{27}{\scriptsize{DRAM:}}} 
          \put(28, 22){\rotatebox{27}{\scriptsize{122.96 GB/sec}}}
        \put(97,101){\scriptsize{DP FMA peak:}}
        \put(95,94){\scriptsize{1327.79 GFLOPS}}
            \end{picture}
       \caption{Parallel}
     \end{subfigure}
      \caption{Roofline plots for (a) serial implementation and (b) threaded implementation.  Threaded implementation here
      uses 28 threads, one for each physical core of the system.  Blue dots ({\color{darkblue}$\bullet$}) are program
      loops and teal squares  ({\color{green}$\blacksquare$})  are loops from Intel MKL routines.  }
           \label{fig:roofline}
   \end{figure}

   The performance of the serial and parallel implementations compared to hardware capabilities are illustrated
   by the roofline plots \cite{williams2009roofline, ilic2014cache} in Figure \ref{fig:roofline}.  The data was collected using Intel Advisor 2019 \cite{intel-advisor} for 
   a tree with 15 levels and $n_c = 16$.  The data was then output using Advisor's \verb|report| option and processed with Tuomas Koskela's \verb|pyAdvisor| tool \cite{pyadvisor-github}. The code is still in the memory-bound region
   for both the serial and parallel implementation, 
   but it is near the transition to compute-bound.  The volume of data copying from children to parent operators required in the merge
   process 
   limits the ability for an implementation to break out of the memory-bound regime.  However, both our serial and parallel implementations are close to the roof.  The serial performance indicates that our implementation is efficient with respect to the hardware's theoretical performance.

\section{Conclusions}
\label{sec:summary}
This paper presented a simple technique for parallelizing the 
two dimensional HPS method for Helmholtz boundary value problems in a shared memory setting with access to parallel linear
algebra.  In the build stage, by far the most computationally expensive stage, we observe a 17.5 times 
speedup over a serial implementation on a desktop computer which is
comparable to a modern super computing node.  
This corresponds to discretizing a problem with over 2 million unknowns and
building the corresponding direct solver in approximately 30 seconds.

While the techniques presented are applied in a shared memory 
setting, they are the foundation for a parallel implementation of the HPS algorithm for the high-frequency Helmholtz equation parallelism on which will be appropriate for upcoming HPC clusters with large memory nodes.  Fully exploiting these machines requires two
level parallelism: using message-passing to divide the computational geometry among several large-memory nodes and shared-memory parallelism, using the techniques presented in this paper, inside the nodes.

The three-dimensional version of the HPS method has much 
larger matrices even close to the leaf level.  Thus  
parallel linear algebra will be utilized earlier in the
build stage and the solve stage will see more benefits
(i.e. larger speedup) from having access to it, as well.

The parallelization technique presented in this 
manuscript can be applied to other tree-based solvers such 
as nested dissection and multifrontal methods.  The matrices 
used to create the calibration data should be modified 
accordingly. For example, the nested dissection method would 
need to use sparse matrices at the leaf level in the build stage.
For other tree-based solvers for dense matrices such as 
 hierarchically semiseparable (HSS) \cite{2007_shiv_sheng,2004_gu_divide}, 
 $\mathcal{H}$-matrix \cite{hackbusch,2003_hackbusch},
hierarchically block separable (HBS) \cite{m2011_1D_survey}
and hierarchical off-diagonal low rank (HODLR) \cite{HOLDR},
it may be necessary to adjust the representative action for a given 
level in the tree.

% We have presented a general and simple technique for optimizing inner and 
% outer threading to acceleration a two-dimensional fast direct solver based 
% on a hierarchical spectral method.  Using the technique, we are able to build 
% the solver in a little over half a minute for a problem with over 2 million unknowns--
% 17.5 times faster than a serial computation--on a modest desktop workstation. 

\section{Acknowledgements}

The authors wish to thank
Total US E\&P for permission to publish.  
The work by A. Gillman is supported by the Alfred P. Sloan foundation and the National Science Foundation (DMS-1522631).  A. Gillman 
and N. Beams are supported in part by a grant from
Total E\&P Research and Technology USA.

\section*{References}
\bibliography{HPSrefs}

\end{document}